\documentclass[aop,seceqn,MSNbibl,citesort,dvips]{arximspdf}
\usepackage{mathrsfs}

\doi{10.1214/10-AOP586}
\volume{39}
\issue{4}
\pubyear{2011}
\firstpage{1468}
\lastpage{1501}

\makeatletter

\newcommand{\implies}{\Longrightarrow}

\newcommand{\del}{\partial}
\newcommand{\delx}{\partial_x}
\newcommand{\delt}{\partial_t}
\newcommand{\lap}{\triangle}
\newcommand{\inv}{^{-1}}

\newcommand{\grad}{\nabla}
\newcommand{\curl}{\grad\times}

\newcommand{\E}{\mathbf{E}}
\newcommand{\as}{\mbox{a.s.}}

\newcommand{\R}{\mathbb{R}}
\newcommand{\Z}{\mathbb{Z}}
\newcommand{\N}{\mathbb{N}}
\newcommand{\T}{\mathbb{T}}

\newcommand{\at}{|}
\newcommand{\varmax}{\vee}
\newcommand{\varmin}{\wedge}
\newcommand{\Chi}[1]{\chi_{_{\{#1\}}}}

\newtheorem{theorem}{Theorem}[section]
\newtheorem{lemma}[theorem]{Lemma}
\newtheorem{proposition}[theorem]{Proposition}

\newproclaim{definition}[theorem]{Definition}
\newproclaim{remark}[theorem]{Remark}

\makeatother

\begin{document}
\begin{frontmatter}

\title{The regularizing effects of resetting in a particle system for the Burgers equation}
\runtitle{The regularizing effects of resetting}

\begin{aug}
\author[A]{\fnms{Gautam} \snm{Iyer}\corref{}\thanksref{t1}\ead[label=e1]{gautam@math.cmu.edu}} and
\author[B]{\fnms{Alexei} \snm{Novikov}\thanksref{t2}\ead[label=e2]{anovikov@math.psu.edu}}
\runauthor{G. Iyer and A. Novikov}
\affiliation{Carnegie Mellon University and Pennsylvania State University}
\address[A]{Department of Mathematical Sciences\\
Carnegie Mellon University\\
Pittsburgh, Pennsylvania 15213\\
USA\\
\printead{e1}}
\address[B]{Department of Mathematics\\
Pennsylvania State University\\
University Park, Pennsylvania 16803\\
USA\\
\printead{e2}}
\end{aug}

\thankstext{t1}{Supported in part by NSF Grant DMS-10-07914 and the
Center for
Nonlinear Analysis.}
\thankstext{t2}{Supported in part by NSF Grant DMS-09-08011.}

\received{\smonth{1} \syear{2010}}

%
\begin{abstract}
We study the dissipation mechanism of a stochastic particle system for
the Burgers equation. The velocity field of the viscous Burgers and
Navier--Stokes equations can be expressed as an expected value of a
stochastic process based on noisy particle trajectories [Constantin and
Iyer \textit{Comm. Pure Appl. Math.} \textbf{3} (2008) 330--345]. In this
paper we study a particle system for the viscous Burgers equations
using a Monte--Carlo version of the above; we consider $N$ copies of
the above stochastic flow, each driven by independent Wiener processes,
and replace the expected value with $\frac{1}{N}$ times the sum over
these copies. A similar construction for the Navier--Stokes equations
was studied by Mattingly and the first author of this paper [Iyer and
Mattingly \textit{Nonlinearity} \textbf{21} (2008) 2537--2553].

Surprisingly, for any finite $N$, the particle system for the Burgers
equations shocks almost surely in finite time. In contrast to the full
expected value, the empirical mean $\frac{1}{N} \sum_1^N$ does not
regularize the system enough to ensure a time global solution. To avoid
these shocks, we consider a resetting procedure, which at first sight
should have no regularizing effect at all. However, we prove that this
procedure prevents the formation of shocks for any $N \geq2$, and
consequently as $N \to\infty$ we get convergence to the solution of
the viscous Burgers equation on long time intervals.
\end{abstract}

%

\end{frontmatter}

\section{Introduction}

The viscous Burgers equation,
%
%
\begin{equation}\label{eqnViscousBurgers}
\delt u + u \,\delx u - \nu\,\delx^2 u = 0,
\end{equation}
has been studied extensively from several different points of view.
Here $\nu> 0$ represents the viscosity, making the equation
dissipative in nature. The inviscid Burgers equation
[equation~(\ref{eqnViscousBurgers}) with $\nu=0$] is studied as the
basic example of a scalar conservation law (see,
e.g.,~\cite{bblEvansPDE,bblDeLellisOttoWestdickenberg}). The Burgers
equation is also linked to the KAM and Aubry--Mather
theories~\cite{bblFathi,bblJauslinKreissMoser}. It is the simplest PDE
that models the Euler and the Navier--Stokes nonlinearity. As such, it
has been extensively studied as the first step in understanding the two
key unresolved issues in fluid mechanics: turbulence and regularity of
the Navier--Stokes equations in three dimensions. In the first category
the objective is to characterize the statistical properties of
turbulence~\cite{bblEKS}. In the second category the objective is to
understand the regularizing mechanism of
dissipation~\cite{bblADV,bblKNS}. This paper falls into the latter
category: we study the regularizing mechanism of a particle system for
the Burgers equations, analogous to the particle system for the
Navier--Stokes equations developed
in~\cite{bblSLNS,bblParticleMethod}.

%

In~\cite{bblSLNS}, a class of second-order nonlinear transport
equations (the Navier--Stokes and viscous Burgers in particular) were
formulated as the average of a stochastic process along noisy particle
trajectories. The formulation for the Navier--Stokes equations
developed in~\cite{bblSLNS} involves recovering the velocity $u$ via
the average of a nonlocal functional of the initial data. For the
viscous Burgers' equation, however, the formulation is simpler.
Explicitly, consider the stochastic flow
%
%
\begin{equation}\label{eqnDXBurgers}
d X_t = u_t (X_t) + \sqrt{2 \nu} \,dW_t
\end{equation}
with initial data $X_0(a) = a$ for all $a \in\R$. Here $W$ denotes a
standard 1D Wiener process. If we require that the velocity $u$ satisfies
%
%
\begin{equation}\label{eqnBurgersUdef}
u_t = \E[ u_0 \circ(X_t^{-1}) ],
\end{equation}
where $\E$ denotes the expected value with respect to the Wiener
measure, then $u$ satisfies\setcounter{footnote}{2}\footnote{This is only valid for spatially
periodic or decay at infinity boundary conditions.} the viscous Burgers
equation~(\ref{eqnViscousBurgers}) and initial data $u_0$. We clarify
that in~(\ref{eqnBurgersUdef}) and subsequently, for any given time $t
\geq0$, $X_t\inv$ denotes the spatial inverse of the diffeomorphism
$X_t$. Namely, we know~\cite{bblKunita}, Theorems 4.5.1, 4.6.5, that
for regular drifts $u$, the stochastic flow $X$ has a modification
which is a stochastic flow of diffeomorphisms of $\R$. Replacing $X$
with this modification if necessary, for any $t \geq0$, we define
$X\inv_t$ to be the inverse of the diffeomorphism $X_t$. That is, for
any $t \geq0$, we have $X_t(X_t\inv(x)) = x$ surely for all $x\in\R$,
and $X_t\inv(X_t(a)) = a$ surely for all $a \in\R$.

Observe that when $\nu= 0$, the system~(\ref{eqnDXBurgers}) and (\ref
{eqnBurgersUdef}) is exactly the
method of characteristics for the inviscid Burgers equation. Indeed
trajectories of the flow $X$ are now characteristics, and equation
(\ref{eqnBurgersUdef}) states that the velocity is transported along
characteristics. Thus, the $\nu> 0$ case could be viewed as a
stochastic generalization of the method of characteristics: we
transport the initial data along noisy characteristics, and then
average with respect to the Wiener measure.

The usual Monte--Carlo method of solving~(\ref{eqnDXBurgers}) and (\ref
{eqnBurgersUdef})
numerically \mbox{\cite
{bblMontecarlo,bblRobertCasella}} is to replace the flow $X_t$ with $N$
different copies $X^{i, N}_t$, each driven by an independent Wiener
process $W^i_t$, and replace the expected value in~(\ref
{eqnBurgersUdef}) by the empirical mean, $\frac{1}{N} \sum_{i=1}^N$.
Explicitly, the system in question becomes
%
%
\begin{eqnarray}
\label{eqndXiN}
dX^{i,N}_t &=& u^N_t(X_t^{i,N}) \,dt + \sqrt{2\nu}
\,dW^i_t,\\
\label{eqnXiN0} X^{i,N}_0(a) &=& a,\\
\label{eqnAn} A^{i,N}_t &=& (X^{i,N}_t)\inv,\\
\label{eqnUn} u^N_t &=& \frac{1}{N} \sum_{i=1}^N u_0 \circ A^{i,N}_t,
\end{eqnarray}
where $u_0$ is the given initial data, $W^i$ a sequence of independent
Wiener processes and $\nu> 0$ the viscosity. As before, for any $t
\geq0$, $(X^{i,N}_t)\inv$ denotes the spatial inverse of $X^{i,N}_t$.
Throughout this paper, with the exception of Section~\ref{sxnShocks},
we impose periodic boundary conditions on the above, and assume the
initial data is periodic with period $1$.

For the Navier--Stokes equations, the particle system in \cite
{bblParticleMethod} involves using a higher-dimensional Wiener process,
and replacing~(\ref{eqnUn}) with the average of vorticity transport and
the Biot--Savart law,
%
%
\begin{eqnarray}
\label{eqnOmegaN} \omega^N_t &=& \E[ ( (\grad X^{i,N}_t)
\omega_0 ) \circ A^{i,N}_t ],\\
\label{eqnBiotSavartN} u^N_t &=& (-\lap)\inv\curl\omega^N_t,
\end{eqnarray}
where $\omega_0 = \curl u_0$ is the initial vorticity. In \cite
{bblParticleMethod}, the authors considered the
system~(\ref{eqndXiN})--(\ref{eqnAn}) and~(\ref{eqnOmegaN}), (\ref
{eqnBiotSavartN}), with spatially periodic boundary conditions,
and proved global existence in two dimensions, local existence in three
dimensions, convergence to the correct limit as $N \to\infty$ and
described the asymptotic behavior for fixed $N$ as $t \to\infty$.

Surprisingly, the techniques of~\cite{bblParticleMethod} \textit{fail}
for the particle system for the Burgers equation [the system
(\ref{eqndXiN})--(\ref{eqnUn})]. Indeed, preliminary numerical
simulations indicate that
the system~(\ref{eqndXiN})--(\ref{eqnUn}) shocks almost surely, in
time \textit
{independent} of $N$. We provide a class of initial data for which we
can prove~(\ref{eqndXiN})--(\ref{eqnUn}) shocks almost surely. We,
however, we are unable
to analytically prove that the shock time is independent of $N$.

One heuristic explanation for the
shock is as follows: this particle system~(\ref{eqndXiN})--(\ref
{eqnUn}) is dissipative
\textit{only} for short time (\cite{bblParticleMethod}, Theorem 5.2).
Once the system~(\ref{eqndXiN})--(\ref{eqnUn}) stops dissipating
energy, the growth from
the nonlinear term should force the system to inherit properties of
the \textit{inviscid} Burgers equation, which shocks if the initial
data is not monotonically nondecreasing.

We remark that the particle system for the Navier--Stokes equations
[the system~(\ref{eqndXiN})--(\ref{eqnAn}) and (\ref
{eqnOmegaN})--(\ref{eqnBiotSavartN})] also dissipates energy \textit{only}
for short time (\cite{bblParticleMethod}, Figure~1 and Theorem 5.2).
However, no dissipation is required to prove 2D global existence for
this system (\cite{bblParticleMethod}, Theorem 3.5). This is because
(\ref{eqnOmegaN}) and~(\ref{eqnBiotSavartN}) are structurally similar to
Euler equations, for which 2D global existence is well known \cite
{bblYudovich} (see also~\cite{bblBealeKatoMajda,bblMajdaBertozzi}). In
contrast, however, the particle system~(\ref{eqndXiN})--(\ref
{eqnUn}) is structurally
similar to the inviscid Burgers equation which is known to shock in
finite time.



The natural approach one would expect to use ``overcoming'' the shocks
in~(\ref{eqndXiN})--(\ref{eqnUn}), would be to continue the system
past shocks as weak
solutions using an analogue of the Rankine--Hugoniot condition (\cite
{bblEvansPDE}, Section 3.4.1), and then prove that as $N \to\infty$,
these weak solutions converge to the smooth solutions of the Burgers
equation. This approach, however, is impossible to use as the
stochastic PDE satisfied by $u^N$ involves second-order terms, for
which the classical techniques (\cite{bblEvansPDE}, Section 3.4.1) will
not work.

While the system~(\ref{eqndXiN})--(\ref{eqnUn}) cannot be continued
past shocks, the
shocks can (surprisingly!) be ``avoided with large probability'' by
resetting the Lagrangian maps. This is the main content of this paper.
Namely, suppose we solve~(\ref{eqndXiN})--(\ref{eqnUn}) for short time
$\delta_t$, and
then replace the initial data with $u^N_{\delta_t}$, and restart the
system~(\ref{eqndXiN})--(\ref{eqnUn}) with this new initial data. Our
main theorem shows
that, if we repeat this procedure often enough, then we can avoid
shocks on an arbitrarily large time interval, with probability
arbitrarily close to $1$.

Explicitly, consider the system
%
%
\begin{eqnarray}
\label{eqndXiNdeltat} dX^{i,N}_{k \delta_t, t} &=& u^N_t(X^{i,N}_{k
\delta_t, t}) \,dt + \sqrt{2\nu} \,dW^i_t,\\
\label{eqndXiN0deltat} X^{i,N}_{k \delta_t, k \delta_t}(a) &=& a,\\
\label{eqnAndeltat} A^{i,N}_{k \delta_t, t} &=& (X^{i,N}_{k\delta_t,
t})\inv,\\
\label{eqnUndeltat} u^N_t &=& \frac{1}{N} \sum_{i=1}^N u^N_{k \delta_t}
\circ A^{i,N}_{k \delta_t, t},
\end{eqnarray}
where $k \in\N$, and $t$ is always assumed to be in the interval $(k
\delta_t, (k+1) \delta_t]$.


If $\delta_t$ is small enough, we show that solutions to this system
exist on arbitrarily large time intervals, with probability arbitrarily
close to one. Once existence with large probability is established, it
is easy to show that as $N \to\infty$ these solutions converge to the
smooth solutions of the viscous Burgers equation.

Before proceeding further, we remark that the fact that the shocks can
be avoided by resetting is doubly unexpected! First, we know that the
inviscid Burgers equations need to be regularized in order for them to
have smooth solutions. The resetting procedure above should morally
provide \textit{no} regularization, as explained in Section~\ref
{sxnGlobalExistence}! Second, the system~(\ref{eqnDXBurgers}) and~(\ref
{eqnBurgersUdef}) is Markovian; if we reset it at regular intervals (as
above), then the new solution obtained will be no different from the
original solution without resetting.

Fortunately, the system~(\ref{eqndXiN})--(\ref{eqnUn}) is not
Markovian, and if we reset
often enough, the generic short time dissipative effect is strong
enough to overcome the nonlinear growth, and with large probability
prevents the formation of shocks. We observed numerically that even large
resetting time $\delta_t$ (i.e., comparable to half the shock time of
the inviscid Burgers system) is enough to ensure that the
system~(\ref{eqndXiNdeltat})--(\ref{eqnUndeltat}) is globally well
posed. With the techniques
in this paper, however, we are only able to prove a global existence
result for~(\ref{eqndXiNdeltat})--(\ref{eqnUndeltat}) when $\delta_t$
is small. The question
for large $\delta_t$ remains open, and cannot be addressed using
techniques in this paper.

Finally, we mention that our technique can be used to show global
existence of the analogue of~(\ref{eqndXiNdeltat})--(\ref
{eqnUndeltat}) for the Navier--Stokes
equations in \textit{two dimensions}. As this is already known \cite
{bblParticleMethod}, \textit{without resetting}, we do not carry out
the details here.

One interesting application would be to the three-dimensional
Navier--Stokes equations. There are numerous results showing global
existence of solutions to the Navier--Stokes equations with small
initial data. One new, interesting question that can be asked in this
framework is the existence of solutions for arbitrary initial data,
which are time global for some small (nonzero) probability. The
McKean--Vlasov-type nonlinearity prevents us from
asking this question for the stochastic Lagrangian formulation for the
Navier--Stokes equations (\cite{bblSLNS}, equations (2.3)--(2.6)). For
the system~(\ref{eqndXiN})--(\ref{eqnAn}) and (\ref
{eqnOmegaN}),~(\ref{eqnBiotSavartN}), the empirical mean $\frac{1}{N}
\sum
_1^N$ provides no regularization, so it is unlikely to expect small
probability time global solutions. However, the repeatedly reset
version of~(\ref{eqndXiN})--(\ref{eqnAn}) and (\ref
{eqnOmegaN}),~(\ref{eqnBiotSavartN}) is free of the McKean--Vlasov
nonlinearity, \textit{and} is dissipative, making it a better candidate
for small probability time global solutions. Unfortunately, there are
obstructions in proving this result directly with the techniques used
here, and we are working on addressing this issue.

The plan of this paper is as follows: in Section \ref
{sxnGlobalExistence} we establish our notational convention, and prove
our main theorem. This proof relies on a few lemmas, the proofs of
which we postpone to Sections~\ref{sxnKeyLemma}, \ref
{sxnFurtherEstimates} and~\ref{sxnCnBound}. In Section~\ref
{sxnShocks}, we provide an example showing that without resetting, the
system~(\ref{eqndXiN})--(\ref{eqnUn}) shocks almost surely. As
mentioned earlier, once
global existence is established the question of convergence as $N \to
\infty$ is easily handled. We conclude the paper by studying this in
Section~\ref{sxnNToInfinity}.

\section{The main theorem and its proof}\label{sxnGlobalExistence}

Throughout this paper, we assume $(\Omega,\Sigma, P)$ is a probability
space and use $\E$ to denote the expected value with respect to the
probability measure $P$.\vadjust{\goodbreak} Let $N \geq2$ be a natural number (which will
be fixed throughout this paper), $\{\mathcal F_t\}_{t\geq0}$ be a
filtration satisfying the usual conditions\footnote{By ``usual
conditions'' (\cite{bblKaratzasShreve}, Definition 2.25) we mean that the
filtration $\{\mathcal F_t\}_{t \geq0}$ is right continuous, and
$\mathcal F_0$ contains all $P$-null sets in $\mathcal F_\infty$.} on
$\Omega$ and $W^1,\ldots, W^N$ be $N$ independent Wiener processes
adapted to the filtration $\{\mathcal F_t\}_{t \geq0}$. We assume
subsequently, without loss of generality, that $\nu= \frac{1}{2}$.

We use $C^k(\T)$, to denote the space of all periodic functions on $\R$
(with period $1$) which have $k$ continuous derivatives. We use $L^p(\T
)$, $H^s(\T)$ to be the Lebesgue $p$-space, and the Sobolev space of
order $s$, respectively, consisting of periodic functions. When writing
norms of functions in these spaces, we drop $\T$. For instance, we use
the notation $\Vert u\Vert_{H^s}$ to denote the $H^s(\T)$ norm of $u$.

We use a calligraphic script to denote the analogous spaces for
processes, on random time intervals. Namely, given $t_0 \geq0$, and a
stopping time $\tau$ such that $\tau\geq t_0$ almost surely, we define
\begin{eqnarray*}
\mathcal C^k([t_0, \tau]; \T) &=& \{ u \mid u \in C^0( [t_0,
\tau]; C^k(\T)) \mbox{ a.s., and } u_{\tau\varmin t} \mbox{ is }
\mathcal F_t \mbox{ adapted} \},\\
\mathcal L^p([t_0, \tau]; \T) &=& \{ u \mid u \in C^0( [t_0,
\tau]; L^p(\T)) \mbox{ a.s., and } u_{\tau\varmin t} \mbox{ is }
\mathcal F_t \mbox{ adapted} \},\\
\mathcal H^s([t_0, \tau]; \T) &=& \{ u \mid u \in C^0( [t_0,
\tau]; H^s(\T)) \mbox{ a.s., and } u_{\tau\varmin t} \mbox{ is }
\mathcal F_t \mbox{ adapted} \},
\end{eqnarray*}
where we use the abbreviation ``a.s.'' for almost surely. For
convenience, if $\tau$ is any stopping time, not necessarily greater
than or equal to $t_0$, we define $\mathcal C^k([t_0, \tau]; \T) =
\mathcal C^k([t_0, \tau\varmax t_0]; \T)$, and similarly for
$\mathcal
H^s, \mathcal L^p$. To avoid confusion with the $C^k(\T)$ norms, we
explicitly use
\[
{\sup_{\omega\in\Omega} \sup_{t_0 \leq t \leq\tau}} \Vert
u_t(\omega)\Vert_{C^{k}}
\]
to denote the $\mathcal C^k([t_0, \tau]; \T)$ norm of $u$.

We clarify, $u \in\mathcal C^k([t_0, \tau]; \T)$ means that there
exists an event $\Omega' \subset\Omega$ with $P(\Omega') = 1$ such
that $\forall\omega\in\Omega'$, $t \in[t_0, \tau(\omega)]$,
$u_t(\omega) \in C^k(\T)$ and $u_t(\omega)$ is continuous in $t$.
Further, $\forall t \geq t_0$, $u_{\tau\varmin t}$ is $\mathcal
F_t$-measurable. In words, $\mathcal C^k( [t_0, \tau]; \T)$ is the set
of all processes which have a $C^k(\T)$ valued, continuous paths
modification and are defined on the random interval $[t_0, \tau]$.

Note that our spaces involve processes which are continuous in time
almost surely, and we are not interested in quantifying any further
regularity with respect to time. When the regularity in time needs to
be quantified, the definition the analogous spaces is not as elementary
(see, e.g.,~\cite{bblKrylovSPDE}).


Our main theorem shows, given any arbitrarily large $T$, we can make our
resetting time $\delta_t$ small enough so that a regular solution
to~(\ref{eqndXiNdeltat})--(\ref{eqnUndeltat}) exists up to time $T$
with probability
arbitrarily close to $1$. In order to formulate our theorem precisely,
we will need to define the notion of solutions to the reset
system~(\ref{eqndXiNdeltat})--(\ref{eqnUndeltat}) with respect to a
stopping time. This is our
next definition.
\begin{definition}\label{dfnSNt0}
Let $t_0 \geq0$, $\tau$ be a \textit{spatially independent} stopping
time such that $\tau\geq t_0$ almost surely, and $u_{t_0}$ be a
$C^1(\T
)$ valued $\mathcal F_{t_0}$-measurable random random variable. Suppose
$u \in\mathcal C^1( [t_0, \tau]; \T)$ is a unique fixed point of the system
%
%
\begin{eqnarray}
\label{eqndXiNt0} X^{i,N}_{t_0, t}(a) &=& a + \int_{t_0}^{\tau\varmin
t} u_s( X^{i,N}_{t_0, s}) \,ds + \int_{t_0}^{\tau\varmin t} d W^i_s,\\
\label{eqnAnt0} A^{i,N}_{t_0, t} &=& ( X^{i,N}_{t_0, t}
)\inv,
\\
\label{eqnUnt0} u_t^{} &=& \frac{1}{N} \sum_{i = 1}^N u_{t_0}^{}
\circ
A^{i,N}_{t_0, t}.
\end{eqnarray}
Then we define
\[
S^{N, \tau}_{t_0, t} u_{t_0} = u_t.
\]
For convenience, we adopt the convention that if $\tau$ is any stopping
time (not necessarily satisfying $\tau\geq t_0$), we define $S^{N,
\tau
}_{t_0, t} u_{t_0} = S^{N, \tau\varmax t_0}_{t_0, t} u_{t_0}$.
\end{definition}
%
%
\begin{remark*}
Note that it is essential to assume $\tau$ does not depend on the
spatial variable, as in this case if the drift $u$ is spatially
regular, then the process $X^{i,N}_{t_0, t}$ admits a modification
which is a stochastic flow of diffeomorphisms. Hence the spatial
inverse $A^{i,N}_{t_0, t} = (X^{i,N}_{t_0, t})\inv$ is well defined. We
will subsequently always assume our stopping times are spatially independent.
\end{remark*}
%
%
\begin{remark*}
In Lemma~\ref{lmaCnBound}, we will show that $S^{N, \tau}_{t_0, \cdot}$
is well defined. Namely, if $k \geq1$, and $u_{t_0} \in C^k(\T)$ is
$\mathcal F_{t_0}$-measurable, then Lemma~\ref{lmaCnBound} shows that
there exists a (deterministic) time $t_1 > t_0$ such that the process
$u$ defined by $u_t = S^{N, t_1}_{t_0, t} u_{t_0}$ belongs to $\mathcal
C^k( [t_0, t_1]; \T)$.
\end{remark*}
%
%
\begin{remark}\label{rmkCompatibility}
Note that we can view the dependence of the operator $S^{N,\tau}_{t_0,
\cdot}$ on the stopping time $\tau$ as a dependence only through the
time interval of definition. Indeed, for any fixed $\delta_t$ and
stopping time $\tau$, $\mathcal C^1([0, \tau]; \T)$ solutions
of~(\ref{eqndXiNt0})--(\ref{eqnUnt0}) are unique up to
indistinguishability.
This follows immediately from a standard argument using
Gronwall's Lemma, and we omit the proof.
Strong uniqueness implies that the operator $S^{N,\tau
}_{t_0, \cdot}$ satisfies a compatibility condition: for $t_0 \geq0$,
consider two stopping times $\tau_1, \tau_2 \geq t_0$ such that $u^1_t
:= S^{N,\tau_1}_{t_0, t} u_{t_0} \in\mathcal C^1([0, \tau_1]; \T)$ and
$u_2 := S^{N,\tau_2}_{t_0, t} u_{t_0} \in\mathcal C^1( [0,\tau_2];
\T
)$, then $u^2 = u^1$ before $\tau_1 \varmin\tau_2$. That is, $u^2$ has
a modification such that for all $t \geq t_0$, $u^2_{t \varmin\tau_1
\varmin\tau_2} = u^1_{t \varmin\tau_1 \varmin\tau_2}$. Thus, when
the time interval of definition is clear, we sometimes omit the
stopping time $\tau$ as a superscript of our operator $S$.
\end{remark}
%

For notational convenience, we omit the first superscript $N$ for the
remainder of this section. Given a (spatially independent) stopping
time $\tau$, a deterministic\vadjust{\goodbreak} starting time $t_0 \geq0$, a $C^1(\T)$
valued $\mathcal F_{t_0}$-measurable initial data $u_{t_0}$, and a
resetting time $\delta_t > 0$ small enough, we can define a $\mathcal
C^1([t_0, \tau]; \T)$ solution of the (stopped) system
(\ref{eqndXiNdeltat})--(\ref{eqnUndeltat}) iteratively by
%
%
\begin{equation}\label{eqnUt}\qquad
\cases{
u_t = u_{t_0}, &\quad when $t = t_0$,\cr
u^{\delta_t}_t = S^\tau_{t_0+ k\delta_t, t} u^{\delta_t}_{t_0 + k
\delta_t}, &\quad whenever $t \in\bigl(t_0 + k\delta_t,
t_0 + (k+1)\delta_t\bigr]$\cr
&\quad for some $k \in\N\cup\{0\}$.}
\end{equation}
%

We are now ready to state our main theorem.
%
%
\begin{theorem}\label{thmGlobalExistence}
Let $N > 1$, $T>0$, $\varepsilon>0$, $s > 6 + \frac{1}{2}$, and
suppose%
\footnote{The theorem and proof remain unchanged if we instead assume
that $u_0$ is a $H^s(\T)$ valued, $\mathcal F_0$-measurable \textit
{bounded} random variable.}
$u_0 \in H^s(\T)$. Then there exists $\delta_T = \delta_T(T,
\varepsilon
, s, \Vert u_0\Vert_{H^{s}})$, independent of $N$, such that for all
$\delta_t < \delta_T$, there exists a spatially independent stopping
time $\tau$ such that $P(\tau> T) > 1 - \varepsilon$ and the process
$u$ defined by~(\ref{eqnUt}) (with $t_0 = 0$) is in the space
$\mathcal
C^6([0, \tau]; \T)$.
\end{theorem}
%
%
\begin{remark}\label{rmkMaximal}
The compatibility condition in Remark~\ref{rmkCompatibility} allows us
to discuss the notion of a maximal stopping time $\tau_{\max}$ for
which the iterative procedure in~(\ref{eqnUt}) is well defined.
Consequently, Theorem~\ref{thmGlobalExistence} will show that this
maximal stopping time $\tau_{\max}$ is in fact at least $T$ with
probability at least $1 - \varepsilon$.
\end{remark}
%

We emphasize that the operator $S_{t_0, t}$ is \textit{not} a smoothing
operator, which, as mentioned earlier, is part of the reason why
Theorem~\ref{thmGlobalExistence} is surprising. We can see $S_{t_0, t}$
is not smoothing from the fact that~(\ref{eqnSPDEu}), the stochastic
partial differential equation (SPDE) satisfied by $u_t = S_{0, t} u_0$
is not dissipative~\cite{bblKrylovSPDE}. One can immediately verify
this as the diffusive term in~(\ref{eqnSPDEu}) does not
necessarily dominate the noise.

Another\vspace*{1pt} (perhaps more intuitive) way of understanding the regularity
properties of $S_{t_0, t}$ is via time splitting. The $S_{t_0}$ can be
time split into two parts: $\bar S^1_{t_0}$, the nonlinear solution
operator associated with the \textit{inviscid} Burgers equations, and
$\bar S^2_{t_0}$ the operator corresponding to resetting. By
considering time split version of~(\ref{eqndXiNdeltat}), one can see
that $\bar S^2_{t_0}$ corresponds exactly to the operator
%
%
\begin{equation}\label{eqnS2barDef}
\bar S^2_{t_0,t_0 + \delta_t} f(x) = \frac{1}{N} \sum_{j=1}^N f
\bigl( x
- ( W^j_{t_0 + \delta_t} - W^j_{t_0} ) \bigr).
\end{equation}

The operator $\bar S^1_{t_0}$ causes growth on the Fourier modes. It is
well known that the damping provided by $\nu\delx^2$, for any $\nu>
0$, is enough to overcome this growth, and this gives us global
existence for the viscous Burgers equations for any strictly positive
viscosity. Thus if the operator $\bar S^2_{t_0}$ provides damping
comparable to $\nu\delx^2$,\vadjust{\goodbreak} then the usual methods can be used to
prove Theorem~\ref{thmGlobalExistence}. However, the operator $\bar
S^2_{t_0}$ provides no damping, as can immediately be checked
from~(\ref{eqnS2barDef}): the operator norm of $\bar S^2_{t_0}$ is
\textit{exactly} $1$ (surely) in all Sobolev and H\"older spaces. This
is the main difficulty in proving Theorem~\ref{thmGlobalExistence}.

We overcome this difficulty by considering the limit $v := \lim
_{\delta
_t \to0} u^{\delta_t}$. It turns out that $v$ satisfies a dissipative
SPDE, and if the initial data is regular enough we obtain convergence
in a strong norm of $u^{\delta_t}$ to $v$. This is the key to the proof
of Theorem~\ref{thmGlobalExistence}, and is formulated below.
\begin{lemma}[(Key lemma)]\label{lmaConvergenceToV}
Let%
\footnote{The proof of this lemma never uses the assumption $N > 1$,
and is valid even for $N = 1$. However, for $N=1$, Lemma \ref
{lmaConvergenceToV} is vacuously true as assumptions (\ref
{eqnUC4bound}) and~(\ref{eqnVC4bound}) will never be satisfied for
nonconstant initial data.}
$N > 1$, $\beta\in\N\cup\{0\}$, $T_0 > t_0 \geq0$, and $\tau$ be a
(spatially independent) stopping time. Let $u_{t_0}$ be a $C^{4+\beta
}(\T)$ valued $\mathcal F_{t_0}$-measurable random variable, and
$u^{\delta_t} \in\mathcal C^{4+\beta}( [t_0, \tau]; \T)$ be defined
by~(\ref{eqnUt}). Let $v \in\mathcal C^{4 + \beta}([t_0, \tau]; \T)$
be the solution of the SPDE
%
%
\begin{equation}\label{eqnSPDEv}
dv_t + v_t \,\delx v_t \,dt - \frac{1}{2} \,\delx^2 v_t \,dt + \frac
{\delx v_t}{N} \sum_{j=1}^N dW^j_t = 0,
\end{equation}
with initial data $v\at_{t=t_0} = u_{t_0}$, and spatially periodic
boundary conditions. Let $\tau_0 = (\tau\varmax t_0) \varmin T_0$, and
$U$ be a constant such that\footnote{Assumptions~(\ref{eqnUC4bound}) and~(\ref{eqnVC4bound}) can be
weakened slightly at the expense of a lengthier, more technical proof.
The weakened assumptions, however, still require more than $\beta$
derivatives. While replacing~(\ref{eqnUC4bound}) and~(\ref{eqnVC4bound})
with a condition involving only $\beta$ derivatives would be of
sufficient interest to warrant a more technical proof, reducing
$4+\beta
$ to $4 + \beta- \varepsilon$ only obscures the heart of the matter.
Since sufficient regularity on our initial data will guarantee (\ref
{eqnUC4bound}) and~(\ref{eqnVC4bound}) anyway, we assume they hold and
avoid unnecessary technicalities.}%
%
%
\begin{eqnarray}
\label{eqnUC4bound}
{\sup_{t_0 \leq t \leq\tau_0} }\Vert u^{\delta_t}_t\Vert_{C^{4+\beta}}
\leq
U \qquad\as,\\
\label{eqnVC4bound}
{\sup_{t_0 \leq t \leq\tau_0} }\Vert v_t\Vert_{C^{4+\beta}} \leq
U \qquad\as,
\end{eqnarray}
and let $w^{\delta_t}_t = u^{\delta_t}_{\tau\varmin t} - v_{\tau
\varmin t}$. Then there exists a constant $C = C(\beta,U,T_0)$,
independent of $N$, $\delta_t$ and $\tau$, such that
%
%
\begin{equation}\label{eqnESTw}
\sup_{t_0 \leq t \leq T_0} \E\Vert w^{\delta_t}_{t}\Vert_{H^\beta}^2
\leq
C \delta_t^{1/2}.
\end{equation}
\end{lemma}
%

Our main interest in this lemma will be for $\beta= 2$, as it will
enable us to obtain a $\mathcal C^1([t_0, T_0]; \T)$ bound on $u$ from
a $\mathcal C^1([t_0, T_0]; \T)$ bound on $v$. A $C^1(\T)$ bound is all
that is needed to continue a solution locally, thus controlling the
$C^1(\T)$ norm of $u$ with large probability, independent of $\delta
_t$, will prove our theorem. Since~(\ref{eqnSPDEv}) is dissipative,
uniform in time bounds of strong norms of $v$ are readily obtained.\vadjust{\goodbreak}
\begin{lemma}\label{lmaH2boundV}
Let $N > 1$, $s \in\N$, $u_0 \in H^s(\T)$. There exists a process $v
\in\mathcal H^s([0, \infty); \T)$ which is a solution to the SPDE
(\ref{eqnSPDEv}) with initial data $u_0$ and periodic boundary
conditions. Further, there exists a constant $V_s =V_s(s, \Vert
u_0\Vert_{H^s})$ such that
%
%
\begin{equation}\label{vBound}
{\sup_{t \geq0} }\Vert v_t\Vert_{H^s} \leq V_s
\end{equation}
almost surely.
\end{lemma}
%

We remark that~(\ref{vBound}) is an almost sure bound on a strong norm
of $v$. The reason we are able to obtain almost sure bounds is because
if we ``multiply by $v$ and integrate by parts'' (or more precisely,
apply It\^o's formula to $\Vert v_t\Vert_{L^2}^2$),\vspace*{1pt} we obtain an equation
with no martingale part! This is carried out in detail in Section~\ref
{sxnFurtherEstimates}.

Lemmas~\ref{lmaConvergenceToV} and~\ref{lmaH2boundV} will now allow
uniform in time control of a strong norm of~$u^{\delta_t}$. The only
remaining ingredient is to obtain a $C^1(\T)$ local existence result,
and guarantee that inequality~(\ref{eqnUC4bound}) is satisfied
uniformly in~$\delta_t$.

\begin{lemma}\label{lmaCnBound}
Suppose $u_{t_0}$ is a $C^1(\T)$ valued $F_{t_0}$-measurable random
variable such that there exists a constant $U^0_1 \geq0$ such that
$\Vert u_{t_0}\Vert_{C^{1}} \leq U^0_1$ almost surely. There exists
$T_0 =
T_0( U^0_1) > t_0$ and a process $u^{\delta_t} \in\mathcal C^1( [t_0,
T_0] ; \T)$ such that $u^{\delta_t}$ is a solution to~(\ref{eqnUt})
with $\tau= T_0$.

If further for some $n \in\N$, $u_{t_0}$ is a $C^n(\T)$ valued
$F_{t_0}$-measurable random variable, and there exists a constant
$U^0_n \geq0$ such that $\Vert u_{t_0}\Vert_{C^{n}} \leq U^0_n$ almost
surely, then $u^{\delta_t}$ is $\mathcal C^n([0,T]; \T)$, and further
there exists a constant $U_n = U_n(U^0_n, n)$, independent of $N$ and
$\delta_t$, such that
%
%
\begin{equation}\label{eqnUdeltatCn}
{\sup_{t_0 \leq t \leq T_0}} \Vert u^{\delta_t}_t\Vert_{C^{n}} \leq
U_n\qquad \as
\end{equation}
for all $\delta_t < T_0$.
\end{lemma}
%
%
\begin{remark*}
The existence time $T_0$ above only depends on a the $C^1(\T)$ norm of
the initial data. However, on the existence interval, any additional
regularity of the initial data is preserved.
\end{remark*}

We are now ready to prove the main theorem. (Lemmas \ref
{lmaConvergenceToV},~\ref{lmaH2boundV} and~\ref{lmaCnBound} will be
proved in Sections~\ref{sxnKeyLemma},~\ref{sxnFurtherEstimates}
and \ref
{sxnCnBound}, resp.)
\begin{pf*}{Proof of Theorem~\ref{thmGlobalExistence}}
Let $\delta_T > 0$ be a small time, to be specified later, and let
$\delta_t \in(0, \delta_T)$ be arbitrary. Given a stopping time
$\tau
$, we define the operator $\mathscr S^{\delta_t, \tau}_{m \delta
_t,t}$ by
\[
\mathscr S^{\delta_t, \tau}_{m \delta_t ,t}=S^\tau_{k \delta_t, t}
\circ S^\tau_{(k-1) \delta_t, k \delta_t} \circ\cdots\circ S^\tau
_{(m+1) \delta_t, (m+2) \delta_t} \circ S^\tau_{m \delta_t, (m+1)
\delta_t},
\]
where $k \in\N$ is such that $k \delta_t < t \leq(k+1)
\delta_t$.\vadjust{\goodbreak}

Let $v_{t}$ be the solution of~(\ref{eqnSPDEv}). By Lemma \ref
{lmaH2boundV} and the Sobolev embedding theorem there is a constant
$V_1$, such that
\[
\sup_{t \geq0} \Vert v_{t}\Vert_{C^{1}} \leq V_1
\]
almost surely. Let $T_0 = T_0(2V_1)$ be the local existence time in
Lemma~\ref{lmaCnBound}; namely, for any initial data $u_0$ with
$\Vert u_0\Vert_{C^{1}} \leq2V_1$, and for any $\delta_t < T_0$, the process
$\mathscr S^{\delta_t, T_0}_{0,\cdot} u_0$ is $\mathcal C^1([0,T_0];
\T
)$. Without loss of generality we can assume that $T_0$ is an integer
multiple of $\delta_t$.

Note that our assumption $u_0 \in H^{13/2+}$, Lemma~\ref{lmaH2boundV}
and the Sobolev embedding theorem imply that assumption (\ref
{eqnVC4bound}) is valid for $\beta= 2$ (in this case, the supremum can
in fact be taken over all $t \in\R$). Similarly, Lemma \ref
{lmaCnBound} guarantees that assumption~(\ref{eqnUC4bound}) is valid
for $\beta= 2$ and all $\delta_t < T_0$. Thus Lemma \ref
{lmaConvergenceToV} can be applied.

Let $\Omega_1$ be the event $\{ \Vert u^{\delta_t}_{T_0}\Vert
_{C^{1}} \leq
2V_1 \}$. Then
\begin{eqnarray*}
P( \Omega_1 ) & \geq& P( \Vert u^{\delta_t}_{T_0} - v_{T_0}
\Vert_{C^{1}}
\leq V_1 )\\
& \geq& P\biggl( \Vert u^{\delta_t}_{T_0} - v_{T_0} \Vert_{H^{2}}
\leq
\frac{V_1}{c_1} \biggr) \qquad\mbox{(Sobolev embedding)}\\
& \geq& 1 - \frac{c_1^2}{V_1^2} \E( \Vert u^{\delta_t}_{T_0}
- v_{T_0} \Vert_{H^{2}}^2 ) \qquad\mbox{(Chebyshev's inequality)}\\
& \geq& 1 - \frac{C \delta_t^{1/2}}{V_1^2} \qquad\mbox{(Lemma \ref
{lmaConvergenceToV})}\\
& \geq& 1 - \frac{C \delta_T^{1/2}}{V_1^2},
\end{eqnarray*}
where the constant $c_1$ above is the constant arising in the Sobolev
embedding theorem. An appropriate choice of $\delta_T$ will make
$P(\Omega_1)$ arbitrarily close to~$1$. We clarify that while our bound
on $P(\Omega_1)$ depends only on $\delta_T$, the event $\Omega_1$
depends on $\delta_t$.

We define a stopping time $\tau_1$ by
\[
\tau_1(\omega) =
\cases{
T_0, &\quad if $\omega\notin\Omega_1$,\cr
2T_0, &\quad if $\omega\in\Omega_1$.}
\]
%
Note that by Remark~\ref{rmkCompatibility}, we have
$\mathscr S^{\delta_t, \tau_1}_{0, t} u_0 =\mathscr S^{\delta_t,
T_0}_{0, t} u_0$ for all $t \in[0, T_0]$. Thus, by the semi-group
property, and the fact that $T_0$ is an integer multiple of $\delta_t$,
%
\[
\mathscr S^{\delta_t, \tau_1}_{0, t} u_0 =
\cases{
\mathscr S^{\delta_t, T_0}_{0, t} u_0, &\quad for $t \in[0, T_0]$,\cr
\mathscr S^{\delta_t, \tau_1}_{T_0, t} \circ\mathscr S^{\delta_t,
T_0}_{0, T_0} u_0, &\quad for $t \in(T_0, 2 T_0]$,}
\]
as long as either side is defined. We claim that the right-hand side
above is well defined and in $\mathcal C^6( [0, \tau_1]; \T)$. We see
this\vspace*{1pt} as follows: first for $t \in[0, T_0]$, this is true by Lemma \ref
{lmaCnBound} and Remark~\ref{rmkCompatibility}. Now, for $\omega\notin
\Omega_1$ and any $t \in[T_0, 2T_0]$, $\mathscr S^{\delta_t,
\tau
_1}_{T_0, t}$ is just the identity operator. Further for almost every
$\omega\in\Omega_1$ we have $\mathscr S^{\delta_t, T_0}_{0,T_0}
u_0(\omega) = u^{\delta_t}_{T_0}(\omega) \in C^6(\T)$ and $\Vert
u^{\delta_t}_{T_0}(\omega)\Vert_{C^{1}} \leq2V_1$. Thus
for almost any $\omega\in\Omega_1$, and for every $t \in[T_0,
2T_0]$, $\mathscr S^{\delta_t, \tau_1}_{T_0, t} u^{\delta
_t}_{T_0}(\omega) \in C^6(\T)$ by Lemma~\ref{lmaCnBound}.

Using Sobolev embedding, Chebyshev's inequality and Lemma \ref
{lmaConvergenceToV} as above, we can find an event $\Omega_2 \subset
\Omega_1$ such that $P(\Omega_2)$ is arbitrarily close to $P(\Omega
_1)$. As before we define a stopping time $\tau_2$ by
\[
\tau_2(\omega) =
\cases{
\tau_1(\omega), &\quad if $\omega\notin\Omega_2$,\cr
3T_0, &\quad if $\omega\in\Omega_2$,}
\]
and the solution $u^{\delta_t}_\cdot= \mathscr S^{\delta_t, \tau
_2}_{0, \cdot} u_0 \in\mathcal C^6( [0, \tau_3]; \T)$. A finite
iteration will complete the proof.
\end{pf*}
%

Finally we address the question of $N \to\infty$. For this purpose, we
re-introduce the superscript of $N$ to indicate the dependence on $N$
of the process considered. Using techniques similar to \cite
{bblParticleMethod}, we show that the solution $v^N$ of (\ref
{eqnSPDEv}) converges to the solution of the viscous Burgers equation as
$N \to\infty$.
\begin{proposition}\label{ppnNToInfinity}
Let $v^{N}$ be the solution of~(\ref{eqnSPDEv}) with initial data
$u_{0}$, and $u^b_{t}$ be the solution of the viscous Burgers
equation~(\ref{eqnViscousBurgers}) with the same initial data. If
$u_{0} \in H^{s}$, $s>\frac{3}{2}$, then for any $T > 0$, there exists
a constant $C = C(T, s, \Vert u_0\Vert_{H^{s}})$ such that
\[
\sup_{t \in[0,T]} \E\Vert u^b_{t}-v^{N}_{t}\Vert_{L^{2}}^2 \leq
\frac{C}{N}.
\]
\end{proposition}
%

We prove Proposition~\ref{ppnNToInfinity} in Section \ref
{sxnNToInfinity}. We conclude by remarking that Proposition \ref
{ppnNToInfinity}, Lemma~\ref{lmaConvergenceToV} and an argument similar
to the proof of Theorem~\ref{thmGlobalExistence} will show that for
small enough $\delta_t$, as $N \to\infty$, $u^{N, \delta_t}$ converges
to the same limit on an event of almost full probability.

\section{Almost sure existence of shocks without resetting}
\label{sxnShocks}

In this section we show that the system~(\ref{eqndXiN})--(\ref
{eqnUn}) develops shocks
almost surely, for any~$N$. The existence of shocks is simpler to prove
if we work with monotone functions on~$\R$, instead of periodic
functions, and thus for this section only, we will work
with~(\ref{eqndXiN})--(\ref{eqnUn}) on $\R$ instead of on $\T$.


Let $\tau$ be a (spatially-independent) stopping time, and we interpret
$\mathcal C^1([0, \tau];\break \R)$ solutions to~(\ref{eqndXiN})--(\ref
{eqnUn}), in the natural
way [analogous to~(\ref{eqndXiNt0})--(\ref{eqnUnt0})]. The main
result of this section
shows that even if we stop ``bad'' realizations of (\ref
{eqndXiN})--(\ref{eqnUn}), we can
never continue solutions past the time $\frac{N}{\Vert\delx u_0 \Vert
_{L^\infty}}$, unless we introduce a regularizing mechanism.\vadjust{\goodbreak}
\begin{proposition}\label{ppnBurgersShock}
Suppose $u_0 \in C^1(\R)$ is a decreasing function, and let $u$ be a
$\mathcal C^1( [0, \tau']; \R)$ a solution of~(\ref{eqndXiN})--(\ref
{eqnUn}) with
initial data $u_0$. Then, almost surely,
%
%
\begin{equation}\label{eqnTauPrimeBound}
\tau' < \frac{N}{\Vert\delx u_0\Vert_{L^\infty}}.
\end{equation}
\end{proposition}
%
%
\begin{remark}\label{rmkShockTime}
The numerically observed shock time, in the periodic case, is
independent of $N$, and it is of the order $\frac{1}{\Vert\delx
u_0\Vert_{L^\infty}}$ with large probability. This indicates our
bound~(\ref{eqnTauPrimeBound}) is far from optimal.
\end{remark}
%
%
\begin{remark}
One can show\footnote{See, for instance, \cite
{bblParticleMethod}, Theorem 4.1, where the analogous result is proved
for the Navier--Stokes equations.} that as $N \to\infty$ the solution
to~(\ref{eqndXiN})--(\ref{eqnUn}) approaches the solution to (\ref
{eqnViscousBurgers}) at
a rate of $\frac{1}{\sqrt{N}}$. However, it is well known that the
solution to~(\ref{eqnViscousBurgers}) is smooth for all time and no
shock develops, provided the initial data is, for instance, $C^1$ and
bounded~\cite{bblEvansPDE}.

The numerics mentioned in Remark~\ref{rmkShockTime}, however, indicate
that no matter how large $N$ is, the system~(\ref{eqndXiN})--(\ref
{eqnUn}) will only be a
good\vspace*{1pt} approximation to the true solution of~(\ref{eqnViscousBurgers})
for short time, in the order of $\frac{1}{\Vert\delx u_0\Vert
_{L^\infty}}$.
\end{remark}
%
%
\begin{remark}
Monotonicity of the initial data $u_0$ is precisely the condition that
constrained us to work on the line instead of on the torus.
Specifically, the assumption $\del_x u_0 (x) < 0$ for arbitrary $x \in
\R$ simplifies the proof of~\ref{ppnBurgersShock} considerably.
Numerics, however, indicate that this monotonicity assumption is
redundant, and~(\ref{eqndXiN})--(\ref{eqnUn}) develops shocks for
arbitrary (periodic)
initial data.
\end{remark}
%
%
\begin{pf*}{Proof of Proposition~\ref{ppnBurgersShock}}
Assume for simplicity, and without loss of generality, that $\Vert
\delx u_0\Vert_{L^\infty} = - \delx u_0(0) = 1$. Let the stopping
time $\tau$ be the
first time $t \leq\tau'$ such that $\delx X^{1,N}_t(0) = 0$. Explicitly,
\[
\tau= \tau' \varmin\inf\{ t \mid \delx X^{1,N}_t(0) = 0 \}.
\]
We will first show that, $\tau\leq N$, almost surely.

Differentiating~(\ref{eqndXiN}) in space gives
%
%
\begin{equation}\label{eqnDGradX}
d ( \delx X^{1,N}_t) = \delx u^N_t \at_{X^{1,N}_t} \delx X^{1,N}_t \,dt,
\end{equation}
for $t < \tau'$, almost surely. Here, our notation $\delx u^N_t\at
_{X^{1,N}_t}$ means
\[
\delx u^N_t\at_{X^{1,N}_t}(x) = \delx u^N_t (X^{1,N}_t(x)).
\]
Differentiating equation~(\ref{eqnUn}) in space, we obtain
%
%
\begin{eqnarray}\label{eqnDiffu}
\delx u^N_t \at_{X^{1,N}_t} &=& \frac{1}{N} \sum_{i=1}^N
\delx
(u_0 \circ A^{i,N}_t) \at_{X^{1,N}_t}\nonumber\\[-8pt]\\[-8pt]
&=& \frac{1}{N} \sum_{i=1}^N \delx u_0 \at_{A^{i,N}_t
\circ X^{1,N}_t} \,\delx A^{i,N}_t \at_{X^{1,N}_t},\nonumber
\end{eqnarray}
for $t < \tau'$ almost surely. Since by the chain rule,
%
%
\begin{equation}\label{eqnGradA1X1}
\delx A^{1,N}_t \at_{X^{1,N}_t} \delx X^{1,N}_t = \delx(
A^{1,N}_t \circ X^{1,N}_t )= 1,
\end{equation}
multiplying~(\ref{eqnDiffu}) by $\delx X^{1,N}_t$ gives
%
%
\begin{eqnarray}\label{eqnTmpGradU1}
&&\delx u^N_t \at_{X^{1,N}_t} \,\delx X^{1,N}_t\nonumber\\[-8pt]\\[-8pt]
&&\qquad = \frac{1}{N} \Biggl[
\delx
u_0 + \sum_{i=2}^N \delx u_0 \at_{A^{i,N}_t \circ X^{1,N}_t} \cdot
\delx A^{i,N}_t\at_{X^{1,N}_t} \,\delx X^{1,N}_t \Biggr]\nonumber
\end{eqnarray}
for $t < \tau'$ almost surely.

Note that for a $C^1$ solution of the system~(\ref{eqndXiN})--(\ref
{eqnUn}), for all $i$,
the flow $X^{i,N}_t\dvtx\R\to\R$ is homotopic to the identity map via
$C^1$ diffeomorphisms of $\R$. The same is true for the inverse inverse
$A^{i,N}_t$, and thus $\delx A^{i,N}_t\at_{X^{1,N}_t} \delx X^{1,N}_t >
0$. Finally, since $u_0$ is assumed to be decreasing, we know that
$\delx u_0 < 0$, and thus equations~(\ref{eqnDGradX}) and (\ref
{eqnTmpGradU1}) yield
\[
\delt\delx X^{1,N}_t(0) < \frac{-1}{N}
\]
for $t < \tau'$ almost surely. This (ordinary) differential inequality,
along with the fact that $\delx X^{1,N}_0 = 1$, necessitates $\tau< N$
almost surely.

Now, by definition of $\tau$, and continuity (in time) of $\delx X^{1,N}_t$,
%
%
\begin{equation}\label{eqnDelxXequals0}
\lim_{t \to\tau^-} \delx X^{1,N}_t(0) = 0
\end{equation}
on the event $\{\tau< \tau'\}$. From~(\ref{eqnTmpGradU1}) and the
chain rule we have
\[
\delx u^N_t \at_{X^{1,N}_t} = \frac{1}{N} \biggl[ \frac{\delx
u_0}{\delx
X^{1,N}_t } + \sum_{i=2}^N \delx u_0 \at_{A^{i,N}_t \circ X^{1,N}_t}
\cdot\delx A^{i,N}_t \at_{X^{1,N}_t} \biggr]
\]
for $t < \tau'$ almost surely. Note that all the terms on the
right-hand side~(\ref{eqnTmpGradU1}) have the same sign. Thus if one of these
terms approaches $-\infty$, then necessarily the entire right-hand side
approaches $-\infty$. Equation~(\ref{eqnDelxXequals0}) immediately
implies the first term approaches $-\infty$ at $x=0$ on the event $\{
\tau< \tau'\}$. Hence, on this event we have
\[
\lim_{t \to\tau^-} \Vert\delx u_t\Vert_{L^\infty} \geq-\lim_{t
\to\tau^-}
\delx u^N_t( X^{1,N}_t(0)) = \infty,
\]
almost surely. Consequently, if $u \in\mathcal C^1([0, \tau']; \R)$,
we must have $P(\tau< \tau') = 0$. Hence $\tau' = \tau< N$ almost surely.
\end{pf*}
%

\section{\texorpdfstring{Proof of the key lemma (Lemma \protect\ref{lmaConvergenceToV})}{Proof of the key lemma (Lemma 2.5)}}
\label{sxnKeyLemma}

In this section we prove convergence of $u^{\delta_t}$ to $v$ as
$\delta
_t \to0$. The basic idea is to show that the velocity in our reset
system~(\ref{eqndXiNdeltat})--(\ref{eqnUndeltat}) satisfies the
limiting SPDE~(\ref{eqnSPDEv})
with a small error which is controlled as $\delta_t \to0$.

By shifting time, we may assume without loss of generality that $t_0 =
0$. Further replacing $\tau$ with $\tau\varmin T_0$ if necessary, we
may assume $\tau= \tau_0 \leq T_0$. Throughout this section, we adopt
the convention that $t_0 = 0$, and $N$, $\beta$, $T_0$, $\tau$, $\tau
_0$, $u_0$, $u^{\delta_t}$, $v$ and $U$ are as in the statement of
Lemma~\ref{lmaConvergenceToV}. We also assume the processes
$X^i_{k\delta_t, \cdot}$, $A^i_{k\delta_t, \cdot}$ are all as in
(\ref{eqndXiNdeltat})--(\ref{eqnUndeltat}), and for notational
convenience, we will omit the
$N$ and $\delta_t$ as superscripts throughout this section.


We need a few lemmas before we can prove Lemma~\ref{lmaConvergenceToV}.
In our first lemma we determine an SPDE satisfied by $u$ on the
interval $(k\delta_t, (k+1) \delta_t]$.
\begin{lemma}\label{lmaSPDEui}
We define the process $u^i$ to be the $i$th summand in (\ref
{eqnUndeltat}). Explicitly,
%
%
\begin{equation}\label{eqnUiDeltatDef}
u^i_t =
\cases{
u_0, &\quad for $t=0$,\cr
u_{k\delta_t} \circ A^i_{k\delta_t, t}, &\quad for $t \in\bigl(\tau\varmin
k\delta_t, \tau\varmin(k+1)\delta_t\bigr]$.}
\end{equation}
Then for all $i \in\{1,\ldots, N\}$, the process $u^i \in\mathcal
C^{4+\beta}([0, \tau]; \T)$, and satisfies the SPDE
%
%
\begin{eqnarray}\label{eqnSPDEui}
&&\Chi{\tau\geq k\delta_t} ( u^i_{\tau\varmin t} - u_{k \delta_t}
) + \int_{\tau\varmin k\delta_t}^{\tau\varmin t} \biggl( u_{s}\,
\delx u^i_{s} - \frac{1}{2} \,\delx^2 u^i_{s}\biggr)\, d{s}
\nonumber\\[-8pt]\\[-8pt]
&&\qquad{}+ \int_{\tau\varmin k\delta_t}^{\tau\varmin t} \delx u^i_{s}
\,dW^i_{s} = 0
\nonumber
\end{eqnarray}
on the interval $t \in[k \delta_t, (k+1)\delta_t]$. Similarly, the
process $u \in\mathcal C^{4+ \beta}([0, \tau]; \T)$, and satisfies
the SPDE
%
%
\begin{eqnarray}\label{eqnSPDEu}
&&u_{\tau\varmin t} - u_{\tau\varmin k\delta_t} + \int_{\tau\varmin
k\delta_t}^{\tau\varmin t} \biggl( u_{s} \delx u_{s} - \frac{1}{2}\,
\delx^2 u_{s}\biggr) \,ds  \nonumber\\[-8pt]\\[-8pt]
&&\qquad{}+ \int_{\tau\varmin k\delta_t}^{\tau\varmin t} \frac{1}{N} \sum
_{j=1}^N \delx u^j_{s} \,dW^j_{s} = 0
\nonumber
\end{eqnarray}
on the interval $t \in[k\delta_t, (k+1)\delta_t]$.
\end{lemma}
%
%
\begin{remark}
A more intuitive, though less precise, way of phrasing the SPDEs (\ref
{eqnSPDEui}) and~(\ref{eqnSPDEu}) would be to say for $t \in(\tau
\varmin k\delta_t, \tau\varmin(k+1)\delta_t]$, $u$, $u^i$ satisfy
the SPDEs
\begin{eqnarray*}
du^i_t + u_t \,\delx u^i_t \,dt - \tfrac{1}{2} \,\delx^2 u^i_t \,dt +
\delx u^i_t \,dW^i_t &=& 0 \qquad\mbox{for all }i \in\{1,\ldots, N\},
\\
du_t + u_t \,\delx u_t \,dt - \frac{1}{2} \,\delx^2 u_t \,dt + \frac
{1}{N} \sum_{j=1}^N \delx u^j_t \,dW^j_t &=& 0
\end{eqnarray*}
with initial data $u^i\at_{t = k\delta_t} = u_{k\delta_t}$ and $u\at_{t
= k\delta_t} = u_{k\delta_t}$.
\end{remark}
\begin{pf*}{Proof of Lemma~\ref{lmaSPDEui}}
From~\cite{bblSLNS,bblParticleMethod} (see also \cite
{bblKrylovRozovski,bblRozovsky}) we know
that when $\tau\equiv\infty$, the process $A^i_{k\delta_t,\cdot}$
satisfies the SPDE
\[
dA^i_{k\delta_t, t} + u_t \,\delx A^i_{k\delta_t, t} \,dt -
\tfrac{1}{2}\,
\delx^2 A^i_{k\delta_t, t} \,dt + \delx A^i_{k\delta_t, t}
\,dW^i_t = 0
\]
on the time interval $(k\delta_t, (k+1)\delta_t]$. Writing down an
integral version of this in the presence of a stopping time, equations
(\ref{eqnSPDEui}) and~(\ref{eqnSPDEu}) follow immediately from (\ref
{eqnUnt0}) and~(\ref{eqnUiDeltatDef}) by a direct application of It\^
o's formula.

To check%
\footnote{The spatial regularity of $u, u^1,\ldots, u^n$ follows
directly from an assumption only on the initial data, and a standard
iteration argument. This is contained in Section~\ref{sxnCnBound}.
However, for Lemma~\ref{lmaSPDEui}, an iteration argument is
unnecessary because of assumption~(\ref{eqnUC4bound}).}
$u, u^1,\ldots, u^N \in\mathcal C^{4+ \beta}([0, \tau]; \T)$, note
that continuity in time is immediate. Further, the spatial regularity
of $u$ has already been assumed in the statement of Lemma \ref
{lmaConvergenceToV}. For $u^1,\ldots, u^N$, note that the $\tau$ and
the noise are spatially independent in~(\ref{eqndXiNdeltat}),
and it immediately shows that each $X^i_{k\delta_t, \cdot}$ (and hence
each $A^i_{k\delta_t, \cdot}$) is as spatially regular as $u$, which in
turn shows that each $u^i \in\mathcal C([0, \tau]; \T)$.
\end{pf*}
%

Now we show that with a small error $u$ satisfies the SPDE (\ref
{eqnSPDEv}) stopped at $\tau$, and obtain bounds on this error. Let
$\E
_{\mathcal F_{k\delta_t}} Y$ denote the conditional expectation of $Y$
given $\mathcal F_{k\delta_t}$. Given any process $f$, and a stopping
time $\tau$, we define the stopped increment $\Delta_k^\tau f$ by
\[
\Delta_k^\tau f = f_{\tau\varmin(k+1)\delta_t} - f_{\tau\varmin
k\delta_t}.
\]
For the (deterministic) process $f_t = t$, we define $\Delta_k^\tau t$ by
\[
\Delta_k^\tau t = \bigl(\tau\varmin(k+1)\delta_t\bigr) - (\tau\varmin
k\delta
_t) =
\cases{
\delta_t, &\quad if $\tau\geq(k+1) \delta_t$,\cr
\tau- k \delta_t, &\quad if $k\delta_t \leq\tau< (k+1) \delta_t$,\cr
0, &\quad if $\tau< k \delta_t$.}
\]
Finally, let $L$ be the (nonlinear) operator defined by
\[
Lu = u\delx u - \tfrac{1}{2} \delx^2 u.
\]

\begin{lemma}\label{lmaEErroru}
Suppose~(\ref{eqnUC4bound}) holds for some $\beta\in\N\cup\{0\}$,
and let $\varepsilon'_k$ be defined by\footnote{In equation (\ref
{eqnEpsilonPrimeDef}), technically, $L u_{k\delta_t}$ is not defined
when $\tau< k \delta_t$. However, in this case, $\Delta_k^\tau t = 0$,
so the value of $L u_{k\delta_t}$ does not matter. We use this
convention subsequently without further mention.}
%
%
\begin{equation}\label{eqnEpsilonPrimeDef}
\varepsilon'_k = \Delta_k^\tau u + Lu_{k \delta_t} \Delta^\tau
_k t +
\delx u_{k\delta_t} \Biggl(\frac{1}{N} \sum_{j=1}^N \Delta_k^\tau
W^j\Biggr).
\end{equation}
Then there exists a constant $C = C(\beta, U, T_0)$ (independent of $N,
k, \delta_t$ and $\tau$) such that for all $\delta_t \leq T_0$ and $k
\leq\frac{T_0}{\delta_t}$ we have
%
%
\begin{eqnarray}
\label{eqnEepsilonPrimeSq} \sup_{x \in\T} \E\vert\delx^\beta
\varepsilon'_k(x)\vert^2 &\leq& C \delta_t^2,\\
\label{eqnEepsilonPrimeGivenFk} \sup_{x \in\T} \E\vert\E
_{\mathcal F_{k\delta_t}} \delx^\beta\varepsilon'_k(x)\vert^2
&\leq& C \delta_t^{3}.
\end{eqnarray}
\end{lemma}
%
%
\begin{remark*}
Since $u$ and all derivatives of $u$ are a priori uniformly bounded
almost surely, the proof of this lemma is straightforward. Without this
a priori bound, we would only obtain similar bounds on $\E\Vert\delx
^\beta\varepsilon'_{k\delta_t}\Vert_{L^2}^2$ and $\E\Vert\E
_{\mathcal F_{k\delta_t}} \delx^\beta\varepsilon'_{k\delta
_t}\Vert_{L^2}^2$, which are
still sufficient for Lemma~\ref{lmaConvergenceToV}.
\end{remark*}
%
%
\begin{pf*}{Proof of Lemma~\ref{lmaEErroru}}
We assume throughout this section that $C$ is a constant only depending
on $U$ and $T_0$ which could change from line to line. Note first that
assumption~(\ref{eqnUC4bound}) and equation~(\ref{eqndXiNt0})
immediately imply that for any $i \in\{1,\ldots, N\}$,
\[
{\sup_{0 \leq t \leq\tau_0}} \Vert\delx X^i_t\Vert_{C^{3+\beta}}
\leq C
\quad\mbox{and}\quad
{\sup_{0 \leq t \leq\tau_0}} \Vert\delx A^i_t\Vert_{C^{3+\beta}}
\leq C
\]
almost surely. Now equation~(\ref{eqnUnt0}) immediately yields the same
bound for $u$, independent of $N$. Thus, making $U$ larger if
necessary, we may assume without loss of generality that (\ref
{eqnUC4bound}) holds for all the processes $u$, $u^1$, $u^2,\ldots, u^n$.

For any $k \in\N\cup\{0\}$, $t \in(k\delta_t, (k+1)\delta_t]$, $n
\leq\beta+ 2$, differentiating~(\ref{eqnSPDEu}) $n$~times gives
\[
\delx^n u_{\tau\varmin t} - \delx^n u_{\tau\varmin k\delta_t} =
-\int
_{\tau\varmin k\delta_t}^{\tau\varmin t} \delx^n Lu_s \,ds - \frac
{1}{N} \sum_{j=1}^N \int_{\tau\varmin k\delta_t}^{\tau\varmin t}
\delx
^{n+1} u^j_s \,dW^j_s,
\]
and hence
\begin{eqnarray*}
&&\E\vert\delx^n u_{\tau\varmin t} - \delx^n u_{\tau\varmin
k\delta_t}\vert^2 \\
&&\qquad\leq2\E\biggl(\int_{\tau\varmin k\delta
_t}^{\tau\varmin t}
\delx^n Lu_s \,ds\biggr)^2\\
&&\qquad\quad{}+ 2\E\Biggl( \frac{1}{N} \sum_{j=1}^N \int_{\tau\varmin k\delta
_t}^{\tau\varmin t} \delx^{n+1} u^j_s \,dW^j_s \Biggr)^2.
\end{eqnarray*}
Note that
\[
\int_{\tau\varmin k\delta_t}^{\tau\varmin t} \delx^{n+1} u^j_s \,dW^j_s
= \int_{k\delta_t}^t \Chi{k\delta_t \leq\tau} \Chi{s \leq\tau}
\delx
^{n+1} u^j_s \,dW^j_s.
\]
Thus using~(\ref{eqnUC4bound}) for both $u$ and $u^i$, for any $t \in
(k \delta_t, (k+1)\delta_t]$ we have
%
%
\begin{equation}\label{eqnSqrtDeltatBoundOnU}
\sup_{x \in\T} \E\vert\delx^n u_{\tau\varmin t}(x) - \delx^n
u_{\tau\varmin k\delta_t}(x)\vert^2 \leq C\Biggl(\delta_t^2 +
\frac{1}{N^2}\sum
_{j=1}^N \delta_t\Biggr) \leq C \delta_t,
\end{equation}
where as usual the constant $C$ may change from line to line, provided
it only depends on $\beta$, $U$ and $T_0$.

Similarly, using~(\ref{eqnSPDEui}) and the above argument we have
%
%
\begin{equation}\label{eqnSqrtDeltatBoundOnUi}
\sup_{x \in\T} \E\Chi{ k\delta_t \leq\tau} \vert\delx^n
u^i_{\tau\varmin t}(x) - \delx^n u_{k\delta_t}(x)\vert^2 \leq C
\delta_t
\end{equation}
for any $t \in(k\delta_t, (k+1)\delta_t]$ and $n \leq2 + \beta$.

Now, from the definition of $\varepsilon'$ and equation (\ref
{eqnSPDEu}) we have
%
%
\begin{eqnarray}\label{eqnEpsilonPrime}
\varepsilon'_k &=&
-\int_{\tau\varmin k\delta_t}^{\tau\varmin(k+1)\delta_t} Lu_s \,ds
- \frac{1}{N} \sum_{j=1}^N \int_{\tau\varmin k\delta_t}^{\tau
\varmin
(k+1) \delta_t} \delx u^j_s \,dW^j_s \nonumber\\
&&{} + Lu_{\tau\varmin k \delta_t} \Delta
_k^\tau
t + \delx u_{\tau\varmin k\delta_t} \Biggl(\frac{1}{N} \sum_{j=1}^N
\Delta_k^\tau W^j\Biggr)\nonumber\\[-8pt]\\[-8pt]
&=& \int_{\tau\varmin k\delta_t}^{\tau\varmin(k+1)\delta_t}
(Lu_{k\delta_t} - Lu_s) \,ds \nonumber\\
&&{} + \frac{1}{N} \sum_{j=1}^N \int
_{\tau\varmin k\delta_t}^{\tau\varmin(k+1)\delta_t} ( \delx
u_{k\delta_t} - \delx u^j_s ) \,dW^j_s\nonumber
\end{eqnarray}
almost surely. For the It\^o integrals in the second term above,
\begin{eqnarray*}
&&\E_{\mathcal F_{k\delta_t}} \biggl( \int_{\tau\varmin k\delta
_t}^{\tau
\varmin(k+1)\delta_t} ( \delx u_{k\delta_t} - \delx u^j_s
)
\,dW^j_s \biggr)\\
&&\qquad=\E_{\mathcal F_{k\delta_t}} \biggl( \int_{ k\delta_t}^{(k+1)\delta_t}
\Chi{\tau\geq k \delta_t} \Chi{s \leq\tau} ( \delx
u_{k\delta_t}
- \delx u^j_{\tau\varmin s} ) \,dW^j_s \biggr)=
0,
\end{eqnarray*}
and hence
\begin{eqnarray*}
\E\vert\E_{\mathcal F_{k\delta_t}} \delx^\beta\varepsilon
'_k\vert^2
&=& \E
\biggl( \delx^\beta\int_{\tau\varmin k\delta_t}^{\tau\varmin
(k+1)\delta_t} (Lu_{k\delta_t} - Lu_s) \,ds\biggr)^2\\
&\leq&\delta_t \int_{k\delta_t}^{(k+1)\delta_t} \E\bigl[ \Chi
{k\delta
_t \leq\tau} \Chi{s \leq\tau} \delx^\beta(Lu_{k\delta_t} -
Lu_{\tau\varmin s})\bigr]^2 \,ds\\
&=& \delta_t \int_{k\delta_t}^{(k+1)\delta_t} \E\bigl[ \Chi{s
\leq\tau
} \delx^\beta(Lu_{\tau\varmin k\delta_t} - Lu_{\tau\varmin
s})\bigr]^2 \,ds\\
&\leq& C \delta_t^3,
\end{eqnarray*}
where the last inequality follows from~(\ref{eqnSqrtDeltatBoundOnU})
with $n=2 + \beta$. This proves~(\ref{eqnEepsilonPrimeGivenFk}).

For~(\ref{eqnEepsilonPrimeSq}), note that the expected value of the
square of the first term in~(\ref{eqnEpsilonPrime}) has already been
bounded by $C \delta_t^3 < C \delta_t^2$. For the second term, the
It\^
o isometry gives
\begin{eqnarray*}
&&\E\Biggl( \frac{1}{N} \sum_{j=1}^N \int_{\tau\varmin k\delta
_t}^{\tau
\varmin(k+1)\delta_t} \delx^\beta( \delx u_{k\delta_t} -
\delx
u^j_s ) \,dW^j_s \Biggr)^2 \\
&&\qquad=\frac{1}{N} \sum_{j=1}^N \int_{k\delta_t}^{(k+1)\delta_t} \E
\bigl[
\Chi{k\delta_t \leq\tau}\Chi{s \leq\tau} \delx^\beta(
\delx
u_{k\delta_t} - \delx u^j_s )\bigr]^2 \,ds,
\end{eqnarray*}
and using~(\ref{eqnSqrtDeltatBoundOnUi}) with $n=1 + \beta$ the proof
is complete.
\end{pf*}
%

We now prove that a time split version of the SPDE~(\ref{eqnSPDEv})
satisfies the same error estimates as in Lemma~\ref{lmaEErroru}.
\begin{lemma}\label{lmaEErrorv}
Suppose~(\ref{eqnVC4bound}) holds for some $\beta\in\N\cup\{0\}$,
and let $\varepsilon''_k$ be defined by
%
%
\begin{equation}\label{eqnEpsilonDoublePrime}
\varepsilon''_k = \Delta_k^\tau v + Lv_{k \delta_t} \Delta
_k^\tau+
\delx v_{k\delta_t} \Biggl(\frac{1}{N} \sum_{j=1}^N \Delta_k^\tau
W^j\Biggr).
\end{equation}
Then bounds~(\ref{eqnEepsilonPrimeSq}) and (\ref
{eqnEepsilonPrimeGivenFk}) hold for $\varepsilon''_k$.
\end{lemma}
%
%
\begin{pf}
First note that
\[
v_{\tau\varmin t} - v_{\tau\varmin k\delta_t} = -\int_{\tau\varmin
k\delta_t}^{\tau\varmin t} Lv_s \,ds - \frac{1}{N} \sum_{j=1}^N
\int
_{\tau\varmin k\delta_t}^{\tau\varmin t} \delx v_s \,dW^j_s
\]
almost surely. Thus for any $n \leq2 + \beta$ and $t \in[k\delta_t,
(k+1)\delta_t]$ using~(\ref{eqnVC4bound}) gives
%
%
\begin{equation}\label{eqnSqrtDeltatBoundOnV}
\sup_{x \in\T} \E\vert\delx^n v_{\tau\varmin t}(x) - \delx^n
v_{\tau\varmin k\delta_t}(x)\vert^2 \leq C \delta_t.
\end{equation}

Similar to the derivation of~(\ref{eqnEpsilonPrime}) we obtain
\[
\varepsilon''_k=
\int_{\tau\varmin k\delta_t}^{\tau\varmin(k+1)\delta_t}
(Lv_{k\delta_t} - Lv_s) \,ds + \frac{1}{N} \sum_{j=1}^N \int
_{\tau\varmin k\delta_t}^{\tau\varmin(k+1)\delta_t} ( \delx
v_{k\delta_t} - \delx v_s ) \,dW^j_s
\]
from definition~(\ref{eqnEpsilonDoublePrime}).
The remainder of the proof is now identical to the proof of Lemma \ref
{lmaEErroru}.\vadjust{\goodbreak}
%
\end{pf}
%

We are now ready to prove Lemma~\ref{lmaConvergenceToV}. We remark that
assumptions~(\ref{eqnUC4bound}) and~(\ref{eqnVC4bound}) are stronger
than necessary. We only need
%
%
\begin{eqnarray}
\label{eqnUVBetaPlus1as} \sup_{0 \leq t \leq\tau} ( \Vert
u_{t}\Vert_{C^{1+\beta}} + \Vert v_{t}\Vert_{C^{1 + \beta}}
) \leq U\qquad
\as,\\
\label{eqnUVBetaPlus2L2} \sup_{x\in\T} \sup_{0 \leq t \leq\tau}
\E
\bigl( \vert\delx^{2+\beta} u_{t}(x)\vert^2 + \vert\delx
^{2+\beta} v_{t}(x)\vert^2 \bigr) \leq U,
\end{eqnarray}
and the bounds on $\varepsilon'$, $\varepsilon''$ provided by
Lemmas~\ref{lmaEErroru} and~\ref{lmaEErrorv} above. The proof we
provide below depends only on these weaker assumptions.
\begin{pf*}{Proof of Lemma~\ref{lmaConvergenceToV}}
Let $\varepsilon_{k} = \varepsilon'_{k} - \varepsilon''_{k}$, where
$\varepsilon'_{k}$, $\varepsilon''_{k}$ are defined by Lemmas \ref
{lmaEErroru} and~\ref{lmaEErrorv}, respectively.
Using~(\ref{eqnEepsilonPrimeSq}), (\ref
{eqnEepsilonPrimeGivenFk}) and the corresponding estimates for
$\varepsilon''_{k}$, we have
%
%
\begin{eqnarray}
\label{eqnEepsilonSq} \sup_{x \in\T} \E\,\delx^\beta\varepsilon_k(x)^2
&\leq& C \delta_t^2,\\
\label{eqnEepsilonGivenFk} \sup_{x \in\T} \E\vert\E_{\mathcal
F_{k\delta_t}} \delx^\beta\varepsilon_k(x)\vert^2 &\leq& C \delta_t^{3}
\end{eqnarray}
for all $k \leq\frac{T_0}{\delta_t}$. As before, we assume $C$ is a
constant that only depends on $\beta$, $U$ and~$T_0$, which may change
from line to line. Now, estimates~(\ref{eqnEepsilonSq}) and~(\ref
{eqnEepsilonGivenFk}) imply
%
%
\begin{eqnarray}
\label{eqnEepsilonSqSob} \E\Vert\varepsilon_{k}\Vert_{H^{\beta}}^2
&\leq&
C \delta_t^2,\\
\label{eqnEepsilonGivenFkSob} \E\Vert\E_{\mathcal F_{k\delta_t}}
\varepsilon_{k}\Vert_{H^{\beta}}^2 &\leq& C \delta_t^{3}.
\end{eqnarray}
For the remainder of the proof we will use the weaker estimates, (\ref
{eqnEepsilonSqSob}) and~(\ref{eqnEepsilonGivenFkSob}).

Now, recall $w_t = w_{\tau\varmin t} =u_{\tau\varmin t} - v_{\tau
\varmin t}$, and we know $w_0 = 0$. Thus
%
%
\begin{eqnarray}\label{eqnDeltakw}
\delx^\beta\Delta_k^\tau w &=& \delx^\beta\Delta_k^\tau
u -
\delx^\beta\Delta_k^\tau v\nonumber\\
&=& -\delx^\beta( L u_{\tau\varmin k \delta_t}
- L v_{\tau\varmin k\delta_t} ) \Delta_k^\tau t\\
&&{} - \delx
^{\beta
+1} w_{k \delta_t} \Biggl(\frac{1}{N} \sum_{j=1}^N \Delta_k^\tau
W^j\Biggr) + \delx^\beta\varepsilon_k.\nonumber
\end{eqnarray}
We first estimate $\E(\delx^\beta\Delta_k^\tau w)^2$ where $k$ is any
integer such that $k \delta_t \leq T_0$.

For this, independence of $W^i$, the mean square of the matringale term
in~(\ref{eqnDeltakw}) is bounded by
%
%
\begin{eqnarray}\label{eqnIndep}\qquad
\E\Biggl[\delx^{\beta+1} w_{k \delta_t} \Biggl(\frac{1}{N}
\sum_{j=1}^N \Delta_k^\tau W^j\Biggr)\Biggr]^2 &=& \E(\delx^{\beta+1}
w_{k\delta_t})^2 \E\Biggl(\frac{1}{N} \sum_{j=1}^N \Delta_k^\tau
W^j\Biggr)^2\nonumber\\[-8pt]\\[-8pt]
&\leq&\frac{\delta_t}{N} \E(\delx^{\beta+1}
w_{k\delta_t})^2.\nonumber
\end{eqnarray}
Next, for the mean square of the first term in~(\ref{eqnDeltakw})
\begin{eqnarray*}
&&\E\bigl(\delx^\beta\bigl( L u_{\tau\varmin k \delta_t}(x) - L
v_{\tau\varmin k\delta_t}(x) \bigr) \bigr)^2 
\\
&&\qquad\leq
C \E[ (\delx^{\beta+2} u_{\tau\varmin k \delta_t})^2 +
(\delx
^{\beta+2} v_{\tau\varmin k \delta_t})^2 ] \\
&&\qquad\quad{}+ C \E[ (\delx^{\beta} (u_{\tau\varmin k \delta_t} \delx
u_{\tau\varmin k \delta_t}))^2+ (\delx^{\beta} (v_{\tau\varmin k
\delta_t} \delx v_{\tau\varmin k \delta_t}))^2 ]
\\
&&\qquad \leq
C \sup_{x \in\T} \max_{0\leq k \leq{T_0}/{\delta_t}} \E
\bigl(
\vert\delx^{2+\beta} u_{\tau\varmin k\delta_t}(x)\vert^2 + \vert
\delx^{2+\beta} v_{\tau\varmin k\delta_t}(x)\vert^2 \bigr) \\
&&\qquad\quad{} + C \max_{0\leq k \leq{T_0}/{\delta_t}} ( \Vert u_{\tau
\varmin k\delta_t}\Vert_{C^{1+\beta}}^4 + \Vert v_{\tau\varmin
k\delta_t}\Vert_{C^{1 + \beta}}^4 ).
\end{eqnarray*}
Hence, using~(\ref{eqnUVBetaPlus1as}) and~(\ref{eqnUVBetaPlus2L2}) we obtain
%
%
\begin{equation}\label{eqnDeltakwTerm1}
\E\bigl(\delx^\beta( L u_{\tau\varmin k \delta_t} - L
v_{\tau
\varmin k\delta_t} ) \Delta_k^\tau t \bigr)^2 \leq\delta
^2_t C.
\end{equation}
By~(\ref{eqnEepsilonSqSob}) the mean square of the last term in (\ref
{eqnDeltakw}) is also bounded by $C \delta_t^{2}$. Thus,
squaring~(\ref{eqnDeltakw}), taking expected values and using Young's
inequality gives
%
%
\begin{equation}\label{eqnEDeltakwsq}
\E(\delx^\beta\Delta_k^\tau w)^2 \leq\frac{3\delta_t}{N} \E
(\delx
^{\beta+1} w_{k\delta_t})^2 + C \delta_t^{2} + 3 \E(\delx
^\beta
\varepsilon_k)^2.
\end{equation}

Now for any $K \leq\frac{T_0}{\delta_t}$,
\begin{eqnarray*}
(\delx^\beta w_{K\delta_t})^2 &=& (\delx^\beta
w_{\tau
\varmin K\delta_t})^2 =
2 \sum_{k=0}^{K-1} \delx^\beta w_{k\delta_t} \,\delx^\beta\Delta
_k^\tau
w + \sum_{k=0}^{K-1} (\delx^\beta\Delta_k^\tau w)^2\\[-3pt]
&=& 2\sum_{k=0}^{K-1} \delx^\beta w_{k\delta_t} \Biggl(
-\delx^\beta(Lu_{\tau\varmin k\delta_t} - Lv_{\tau\varmin
k\delta_t}) \Delta_k^\tau t \\[-3pt]
&&\hspace*{62.3pt}{} - \delx^{\beta+ 1} \Delta_k^\tau w \Biggl(\frac{1}{N} \sum_{j=1}^N
\Delta_k^\tau W^j \Biggr) + \delx^\beta\varepsilon_k \Biggr)
\\[-3pt]
&&{} +\sum_{k=0}^{K-1} (\delx^\beta\Delta_k^\tau w)^2.
\end{eqnarray*}
Taking expected values, integrating in space using (\ref
{eqnEDeltakwsq}) and~(\ref{eqnEepsilonSqSob}) gives
%
%
\begin{eqnarray}\label{eqnEwKdeltatL2sq}\quad
\E\Vert\delx^\beta w_{K\delta_t}\Vert_{L^2}^2 &\leq&
-2 \delta_t \sum_{k=0}^{K-1} \E\int_\T\delx^\beta w_{k\delta_t}
\delx
^\beta(
u_{\tau\varmin k\delta_t} \delx u_{\tau\varmin k\delta_t}\nonumber\\[-3pt]
&&\hspace*{109.7pt}{}
- v_{\tau\varmin k\delta_t} \delx v_{\tau\varmin
k\delta_t}) \,dx \\[-3pt]
&&{} + 2 \sum_{k=0}^{K-1} \E\int_\T\delx^\beta w_{k\delta_t}\, \delx
^\beta
\varepsilon_k \,dx \nonumber\\[-3pt]
&&{} - \biggl(1 - \frac{1}{N}\biggr) \delta_t \sum_{k=0}^{K-1} \E\int
_\T
(\delx^{\beta+1} w_{k\delta_t})^2 \,dx + C K \delta
_t^{2}.\nonumber
%
\end{eqnarray}
For the first term on the right-hand side of inequality (\ref
{eqnEwKdeltatL2sq}) note
%
%
\begin{eqnarray}\label{eqnEwKdeltatL2sqTerm1}
&&\delx^\beta w_{k\delta_t} \delx^\beta( u_{\tau\varmin
k\delta_t}\,
\delx u_{\tau\varmin k\delta_t} - v_{\tau\varmin k\delta_t} \,\delx
v_{\tau\varmin k\delta_t}) \nonumber\\[-8pt]\\[-8pt]
&&\qquad=\delx^\beta w_{k\delta_t} \,\delx^\beta( w_{k\delta_t} \,\delx
u_{\tau\varmin k\delta_t} - v_{\tau\varmin k\delta_t} \,\delx
w_{k\delta
_t}).
\nonumber
\end{eqnarray}
Observe that the mass (spatial mean) of solutions to (\ref
{eqnSPDEv}) is constant in time. The same is true for solutions
to~(\ref{eqnSPDEu}). Thus, for all $t \leq T_0$, $\int_\T u_{\tau
\varmin t}
\,dx = \int_\T u_0 \,dx = \int_\T v_{\tau\varmin t} \,dx$, and hence
$\int_\T w_t \,dx = 0$. Thus integrating~(\ref{eqnEwKdeltatL2sqTerm1})
in space and using the Poincar\'e inequality, the term involving $u$
above is bounded by
\[
\biggl\vert\int_\T\delx^\beta w_{k\delta_t} \,\delx^\beta( w_{k
\delta_t} \delx u_{\tau\varmin k \delta_t} ) \,dx \biggr\vert
\leq C \Vert\delx^\beta w_{k\delta_t}\Vert_{L^2}^2 \Vert u_{\tau
\varmin k\delta_t}\Vert_{C^{\beta+1}}.
\]
For the term involving $v$ in~(\ref{eqnEwKdeltatL2sqTerm1}), when all
the derivatives fall on $w$ we have
\[
\delx^\beta w_{k\delta_t} v_{\tau\varmin k \delta_t} \,\delx^{\beta+1}
w_{k \delta_t} = \tfrac{1}{2} v_{\tau\varmin k\delta_t} \,\delx(
\delx^\beta w_{k\delta_t})^2,
\]
and if we integrate by parts, we can avoid the extra derivative on $w$. Thus
\[
\biggl\vert\int_\T\delx^\beta w_{k\delta_t} \,\delx^\beta(
v_{\tau\varmin k \delta_t} \,\delx w_{k \delta_t} ) \,dx \biggr\vert
\leq C \Vert\delx^\beta w_{k\delta_t}\Vert_{L^2}^2 \Vert v_{\tau
\varmin k\delta_t}\Vert_{C^{\beta}}.
\]
Thus using~(\ref{eqnUVBetaPlus1as}) and the above estimates, the first
term on the right-hand side of~(\ref{eqnEwKdeltatL2sq}) is bounded by
%
%
\begin{eqnarray}\label{eqnTerm1Ewk}
&&-2\delta_t \E\int_\T\delx^\beta w_{k\delta_t} \,\delx^\beta(
u_{\tau\varmin k\delta_t} \,\delx u_{\tau\varmin k\delta_t} - v_{\tau
\varmin k\delta_t} \,\delx v_{\tau\varmin k\delta_t} )
\,dx\nonumber\\[-8pt]\\[-8pt]
&&\qquad\leq
C \delta_t \E\Vert\delx^\beta w_{k\delta_t}\Vert_{L^2}^2.\nonumber
\end{eqnarray}

For the second term in~(\ref{eqnEwKdeltatL2sq}), we know $w_{k\delta
_t}$ is $\mathcal F_{k\delta_t}$-measurable. Thus using (\ref
{eqnEepsilonGivenFkSob}) and the Cauchy--Schwarz inequality we obtain
%
%
\begin{equation}\label{eqnTerm2Ewk}
\int_\T\vert\E\delx^\beta w_{k\delta_t} \,\delx^\beta\varepsilon
_k\vert\,
dx = \int_\T\vert\E( \delx^\beta w_{k\delta_t} \,\E
_{\mathcal F_{k\delta_t}} \delx^\beta\varepsilon_k )\vert\,
dx \leq C
\delta_t^{3/2}.
\end{equation}

The third term on the right-hand side of~(\ref{eqnEwKdeltatL2sq}) is always
nonpositive, and can be ignored. Thus, recalling $K \leq\frac
{T_0}{\delta_t}$, and using~(\ref{eqnTerm1Ewk}) and (\ref
{eqnTerm2Ewk}) in~(\ref{eqnEwKdeltatL2sq}) we have
\[
\E\Vert\delx^\beta w_{K\delta_t}\Vert_{L^2}^2 \leq C \delta
_t^{1/2} + C
\sum_{k=0}^{K-1} \E\Vert\delx^\beta w_{k\delta_t}\Vert
_{L^2}^2\delta_t.
\]

The remainder of the proof is an elementary discrete Gronwall argument. Let
\[
y_K = C \delta_t^{1/2} + C \sum_{k=0}^{K-1} \E\Vert\delx^\beta
w_{k\delta_t}\Vert_{L^2}^2\delta_t.
\]
Then
\[
y_{k+1} - y_k = C \delta_t \E\Vert\delx^\beta w_{k\delta_t}\Vert_{L^2}^2
\leq C \delta_t y_k
\]
and hence
\[
y_{k+1} \leq( 1 + C \delta_t ) y_k.
\]
Iterating, and using $y_0 = C \delta_t^{1/2}$ gives
\[
y_k \leq(1 + C \delta_t )^k C \delta_t^{1/2}.
\]
Since $k \leq\frac{T_0}{\delta_t}$ this gives

%
%
\begin{equation}\label{eqnYkFinal}
\max_{k \leq{T_0}/{\delta_t}} y_k \leq C \delta_t^{1/2} \sup
_{\delta_t' > 0} ( 1 + C \delta_t') ^ {T_0 / \delta_t'}
\leq C \delta_t^{1/2} e^{C T_0}.
\end{equation}
This proves~(\ref{eqnESTw})\vspace*{1pt} for all times $t$ which are an integer
multiple of $\delta_t$. Since for any $x \in\T$, and $k \leq\frac
{T_0}{\delta_t}$ we elementarily have
\[
\sup_{k \delta_t \leq t \leq(k+1) \delta_t} \E\vert\delx^\beta
v_{\tau\varmin t}(x) - \delx^\beta v_{\tau\varmin k\delta
_t}(x)\vert^2
\leq C \delta_t
\]
and
\[
\sup_{k \delta_t \leq t \leq(k+1) \delta_t} \E\vert\delx^\beta
u_{\tau\varmin t}(x) - \delx^\beta u_{\tau\varmin k\delta
_t}(x)\vert^2
\leq C \delta_t
\]
completing the proof.
\end{pf*}
%

\section{\texorpdfstring{Proof of Lemma \protect\ref{lmaH2boundV}}{Proof of Lemma 2.6}}
\label{sxnFurtherEstimates}
In this section we establish uniform in time bounds for the solution
of~(\ref{eqnSPDEv}) and prove as in Lemma~\ref{lmaH2boundV}. We do this
via the following two lemmas:
\begin{lemma}\label{lmaShortTimeBoundsOnV}
Let $u_0 \in C^\infty(\T)$, $T > 0$, and suppose $v \in\mathcal
C^\infty([0, T]; \T)$ is a solution to~(\ref{eqnSPDEv}) with initial
data $u_0$ and periodic boundary conditions. Then for any $s \in\Z^+$,
there exists a constant $V_s = V_s(s, T, \Vert u_0\Vert_{H^{s}})$ such that
\[
{\sup_{0 \leq t \leq T}} \Vert v_t\Vert_{H^{s}} \leq V_s
\]
almost surely.
\end{lemma}
%
%
\begin{lemma}\label{lmaLongTimeBoundsOnV}
Let $u_0 \in C^\infty(\T)$, and suppose $v \in\mathcal C^\infty
([0,\infty), \T)$ is a solution to~(\ref{eqnSPDEv}) with initial data
$u_0$ and periodic boundary conditions. Then for any $s \in\Z^+$, $T >
0$, there exists a constant $V_s = V_s(s, T, \Vert u_0\Vert_{L^2})$
such that
%
%
\begin{equation}\label{eqnLongTimeHsBoundV}
\sup_{t \geq T} \Vert v_t\Vert_{H^{s}} \leq V_{s}
\end{equation}
almost surely.
\end{lemma}

We draw attention to the fact that a priori bounds are almost sure!
Indeed, applying It\^o's formula to $\Vert v_t\Vert_{L^2}^2$ immediately
yields an equation with no martingale part [see (\ref
{eqnL2normV}) below].

Given Lemmas~\ref{lmaShortTimeBoundsOnV} and \ref
{lmaLongTimeBoundsOnV}, the proof of Lemma~\ref{lmaH2boundV} is now immediate.
\begin{pf*}{Proof of Lemma~\ref{lmaH2boundV}}
Given the almost sure a priori bounds in Lem\-mas~\ref
{lmaShortTimeBoundsOnV} and~\ref{lmaLongTimeBoundsOnV}, existence of
solutions to~(\ref{eqnSPDEv}) follows via standard methods. The time
global bound~(\ref{vBound}) is also an immediate consequence of
Lemmas~\ref{lmaShortTimeBoundsOnV} and~\ref{lmaLongTimeBoundsOnV}.
\end{pf*}
%

We devote the remainder of this section to proving Lemmas \ref
{lmaShortTimeBoundsOnV} and~\ref{lmaLongTimeBoundsOnV}.
\begin{pf*}{Proof of Lemma~\ref{lmaShortTimeBoundsOnV}}
We prove Lemma~\ref{lmaShortTimeBoundsOnV} via energy estimates. First
note that It\^o's formula and~(\ref{eqnSPDEv}) give
\begin{eqnarray*}
d(v_t)^2 &=& 2v_t \,dv_t + \frac{1}{N^2}\sum_{j=1}^N (\delx v_t)^2
\,dt\\
&=& -2 v_t^2 \,\delx v_t \,dt + v_t \,\delx^2 v_t \,dt -2 \frac{v_t\,
\delx
v_t}{N} \sum_{j=1}^N dW^j_t\\
&&{} + \frac{1}{N} (\delx v_t )^2 \,dt.
\end{eqnarray*}
Integrating in space, and using $\int_\T v_t \,\delx v_t \,dx = 0 =
\int
_\T v_t^2 \,\delx v_t \,dx$ gives
%
%
\begin{equation}\label{eqnL2normV}
\del_t \Vert v_t\Vert_{L^2}^2 = -\biggl(1 - \frac{1}{N} \biggr)
\Vert\delx v_t\Vert_{L^2}^2 \quad\implies\quad\Vert v_t\Vert
_{L^2} \leq\Vert u_0\Vert_{L^2}
\end{equation}
almost surely.

A similar calculation shows $\Vert v_t\Vert_{L^p} \leq\Vert u_0\Vert
_{L^p} $ for
all $p \geq2$, and hence\footnote{This can alternately be shown using
a version of the maximum principle~\cite{bblKrylovMaximumPrinciple}.}
$\Vert v_t\Vert_{L^\infty} \leq\Vert u_0\Vert_{L^\infty}$. Recall
$s \geq1$ by assumption,
and so the Sobolev embedding theorem shows $\Vert u_0\Vert_{L^\infty}
\leq c
\Vert u_0\Vert_{H^{s}}$ for some absolute constant $c$.

Now, differentiating~(\ref{eqnSPDEv}) with respect to $x$ and applying
It\^o's formula to $(\del_{x} v_t)^{2}$ we obtain
\begin{eqnarray*}
d (\del_{x} v_t)^{2} &=& 2 \,\del_{x} v_t \,d (\del_{x} v_t) + \frac
{1}{N}\vert\del^{2}_x v_t\vert^{2} \,dt\\
&=& - 2 \,\del_{x} v_t \Biggl( \del_{x}(v_t \delx v_t)
\,dt -
\frac{1}{2} \,\delx^3 v_t \,dt + \frac{\del^{2}_x v_t}{N} \sum_{j=1}^N
dW^j_t \Biggr)\\
&&{}+\frac{1}{N}\vert\del^{2}_x v_{t}\vert^{2} \,dt.
\end{eqnarray*}
Integrating with respect to $x$ on $[0,1]$, and noting that $\int_{\T}
\del_{x} v_t \,\del^{2}_x \,v_t \,dx = 0$, gives
%
%
\begin{eqnarray}\label{eqnIneq}
&&d \Vert\del_{x} v_{t}\Vert_{L^2}^{2} = - \biggl(1 - \frac{1}{N}\biggr)
\Vert\del_{x}^{2} v_{t}\Vert_{L^2}^{2} \,dt + \biggl(2 \int_0^1
\delx^2
v_t (v_t \,\delx v_t) \,dx\biggr) \,dt\nonumber\\
&&\quad
\implies\quad\del_t \Vert\del_{x} v_{t}\Vert_{L^2}^{2} \leq- \frac{1}{4}
\Vert\del_{x}^{2} v_{t}\Vert_{L^2}^{2} + 8 \Vert v_{t} \del_{x}
v_{t}\Vert_{L^2}^{2}\nonumber\\[-8pt]\\[-8pt]
&&\qquad\hspace*{2.6pt}\phantom{\implies d \Vert\del_{x} v_{t}\Vert_{L^2}^{2}}
\leq- \frac{1}{4} \Vert\del_{x}^{2} v_{t}\Vert_{L^2}^{2} + 8
\Vert v_{t}\Vert_{L^\infty}^{2} \Vert\del_{x} v_{t}\Vert
_{L^2}^{2}\nonumber\\
&&\qquad\hspace*{2.6pt}\phantom{\implies d \Vert\del_{x} v_{t}\Vert_{L^2}^{2}}
\leq- \frac{1}{4} \Vert\del_{x}^{2} v_{t}\Vert_{L^2}^{2} + 8
\Vert u_{0}\Vert_{L^\infty}^{2} \Vert\del_{x} v_{t}\Vert
_{L^2}^{2},\nonumber
\end{eqnarray}
almost surely.
Thus,~(\ref{eqnL2normV}),~(\ref{eqnIneq}) and Gronwall's inequality gives
%
%
\begin{equation}\label{eqnH1boundV}
\Vert v_{t}\Vert_{H^{1}} \leq C_{1} e^{c_{0}t}
\quad\mbox{and}\quad
\int_{0}^{t}\Vert v_{t'}\Vert_{H^{2}}^{2} \,dt' \leq C_{1} e^{c_{0}t},
\end{equation}
almost surely, for some constants $C_{1}=C_{1}(\Vert u_{0}\Vert_{H^{1}})$
and $c_{0}=c_{0}(\Vert u_{0}\Vert_{L^{\infty}})$.

For the remainder of this proof we adopt the convention that $c$, $C$
denote absolute constants, $C_s = C_s(s, \Vert u_0\Vert_{H^{s}})$ denotes
a constant depending only on $s$, $\Vert u_0\Vert_{H^{s}}$ and $c_0$
denotes a constant depending only on $\Vert u_0\Vert_{L^\infty}$. The
exact value
of these constants are immaterial, and we will allow them to change
from line to line.

Similar to~(\ref{eqnIneq}), differentiating~(\ref{eqnSPDEv}) twice
with respect to $x$, applying It\^o's formula to $(\del^{2}_{x}
v_t)^{2}$, integrating in space, noting $\int_\T\delx^2 v_t \,\delx^3
v_t \,dx = 0$ and using H\"older's inequality gives
\begin{eqnarray*}
\del_t \Vert\del^{2}_{x} v_{t}\Vert_{L^2}^{2} &\leq&- \biggl( 1
-\frac
{1}{N}\biggr)\Vert\del^{3}_{x} v_{t}\Vert_{L^2}^{2} +2 \Vert\del
^{3}_{x} v_{t}\Vert_{L^2} \Vert\del_{x} (v_{t} \del_{x} v_{t}
) \Vert_{L^2}\\
& \leq&
- c \Vert\del_{x}^{3} v_{t}\Vert_{L^2}^{2}
+ C ( \Vert\del_{x} v_{t}\Vert_{L^\infty}^{2} + \Vert
v_{t}\Vert_{L^\infty}^{2} )\Vert\del^{2}_{x} v_{t}\Vert
_{L^2}^{2}\\
&\leq&- c \Vert\del_{x}^{3} v_{t}\Vert_{L^2}^{2} + C \Vert
v_{t}\Vert_{H^{2}}^{2} \Vert\del^{2}_{x} v_{t}\Vert_{L^2}^{2},
\end{eqnarray*}
almost surely, where the last inequality is obtained by the Sobolev
embedding theorem.
Using~(\ref{eqnH1boundV}), this gives
\[
\Vert v_{t}\Vert_{H^{2}} \leq C_{2} e^{\int_{0}^{t}\Vert v_{t}\Vert
_{H^{2}}^{2}
\,dt} \leq C_{2} e^{C_{1}e^{c_{0}t}}
\quad\mbox{and}\quad
\int_{0}^{t}\Vert v_{t'}\Vert_{H^{3}}^{2} \,dt' \leq C_{2}
e^{C_{1}e^{c_{0}t}},
\]
almost surely.
Proceeding inductively, suppose we know
%
%
\begin{equation}\label{eqnHsBoundV}
\cases{
\Vert v_{t}\Vert_{H^{s}} \leq
C_{s} \exp\bigl( C_{s-1}\exp\bigl( C_{s-2} \cdots\exp
(c_{0}t) \cdots\bigr) \bigr),\cr
\displaystyle \int_{0}^{t}\Vert v_{t'}\Vert_{H^{s+1}}^{2} \,dt' \leq C_{s} \exp
\bigl(
C_{s-1}\exp\bigl( C_{s-2} \cdots\exp(c_{0}t) \cdots
\bigr) \bigr),
}
\end{equation}
holds almost surely for some $s \in\Z^+$. Differentiating (\ref
{eqnSPDEv}) $s+1$ times with respect to $x$, applying It\^o's formula
for $(\del^{s+1}_{x} v_t)^{2}$ and integrating in space we obtain
\[
d \Vert\del^{s+1}_{x} v_{t}\Vert_{L^2}^{2} = - \biggl( 1 -\frac
{1}{N}
\biggr)\Vert\del^{s+2}_{x} v_{t}\Vert_{L^2}^{2} \,dt + 2 \Vert\del
^{s+2}_{x} v_{t}\Vert_{L^2} \Vert\del^{s}_{x} (v_{t} \del_{x}
v_{t} ) \Vert_{L^2} \,dt
\]
since $\int_\T\delx^{s+1} v_t \,\delx^{s+2} v_t \,dx = 0$. Thus
\begin{eqnarray*}
\del_t \Vert\del^{s+1}_{x} v_{t}\Vert_{L^2}^{2} &\leq&- c \Vert
\del_{x}^{s+2} v_{t}\Vert_{L^2}^{2} \\
&&{}+
C ( \Vert\del^{s}_{x} v_{t}\Vert_{L^\infty}^{2} \Vert\del
_{x} v_{t}\Vert_{L^2}^{2}+ \cdots\\
&&\hspace*{24.6pt}{}+ \Vert\del_{x} v_{t}\Vert_{L^\infty}^{2} \Vert\del^{s}_{x}
v_{t}\Vert_{L^2}^{2} + \Vert v_{t}\Vert_{L^\infty}^{2} \Vert\del
^{s+1}_{x} v_{t}\Vert_{L^2}^{2} )
\\
& \leq&- c \Vert\del_{x}^{s+2} v_{t}\Vert_{L^2}^{2}\\
&&{} + C (\Vert
\del^{s}_{x} v_{t}\Vert_{L^\infty}^{2} + \cdots \\
&&\hspace*{24.6pt}{} +
\Vert\del_{x} v_{t}\Vert_{L^\infty}^{2}+ \Vert v_{t}\Vert_{L^\infty}^{2} )\Vert\del^{s+1}_{x}
v_{t}\Vert_{L^2}^{2}
\\
&\leq&- c \Vert\del_{x}^{s+2} v_{t}\Vert_{L^2}^{2} + C \Vert
v_{t}\Vert_{H^{s+1}}^{2} \Vert\del_{x}^{s+1} v_{t}\Vert_{L^2}^{2},
\end{eqnarray*}
almost surely. Thus by Gronwall's lemma
\begin{eqnarray*}
\Vert v_{t}\Vert_{H^{s+1}} &\leq& C_{s+1} \exp\biggl(\int
_{0}^{t}\Vert v_{t'}\Vert_{H^{s+1}}^{2} \,dt'\biggr)\\
&\leq& C_{s+1} \exp\bigl( C_{s}\exp( C_{s-1} \cdots\exp
(c_{0}t) \cdots) \bigr)
\end{eqnarray*}
almost surely. Further
\[
\int_{0}^{t}\Vert v_{t'}\Vert_{H^{s+2}}^{2} \,dt' \leq C_{s+1} \exp
\bigl(
C_{s}\exp( C_{s-1} \cdots\exp(c_{0}t) \cdots
) \bigr),
\]
almost surely, completing the inductive step. By induction, (\ref
{eqnHsBoundV}) holds for all $s \in\Z^+$ completing the proof.
\end{pf*}
%
%
\begin{pf*}{Proof of Lemma~\ref{lmaLongTimeBoundsOnV}}
We prove Lemma~\ref{lmaLongTimeBoundsOnV} via a bootstrapping argument
in Fourier space. To fix notation, for $n\! \in\!\Z$, we use $\hat f(n)
\!=\!
\int_\T e^{-2\pi i n x} f(x) \,dx$ to denote the $n$th Fourier
coefficient of $f$.

On Fourier coefficients, using $u \,\delx u = \frac{1}{2} \,\delx u^2$,
equation~(\ref{eqnSPDEv}) reduces to
%
%
\begin{eqnarray}\label{eqnSdeFourier}
&&
d \hat{v}_t (n) + \frac{2 \pi i n}{N} \hat{v}_t (n) \sum_{j=1}^N d
W^j_t \nonumber\\[-8pt]\\[-8pt]
&&\qquad{}+ 2 \pi^2 n^2 \hat{v}_t (n) \,dt + \pi i n\sum_{m \in\Z}
\hat
v_t(n-m) \hat v_t(m) \,dt = 0\nonumber
\end{eqnarray}
for every $n \in\Z$.

By It\^o's formula applied to~(\ref{eqnSdeFourier})
%
%
\begin{eqnarray}\label{eqndVHatNSq1}
d \vert\hat{v}_t (n)\vert^{2} &=& \overline{ \hat v_t(n) } \,d \hat{v}_t
(n) \,dt+ \hat{v}_t (n) d \,\overline{ \hat v_t(n) }+\frac{4 \pi^{2}
n^{2}}{N} \vert\hat{v}_t (n)\vert^{2} \,dt\nonumber\\
& = & - 4 \pi^{2} n^{2} \biggl( 1 -\frac{1}{N}\biggr) \vert\hat{v}_t
(n)\vert^{2} \,dt \\
&&{} + \pi i n \bigl( \hat{v}_t (n) \overline
{{B}_{t}(n)} -
\overline{ \hat v_t(n) }B_{t}(n) \bigr) \,dt,\nonumber
\end{eqnarray}
where $\overline{\hat v_t (n) }$ denotes the complex conjugate of
$\hat
v_t(n)$, and
\[
B_{t}(n) = \sum_{m \in\Z} \hat v_t(n-m) \hat v_t(m),
\]
is the nonlinear Fourier coupling in~(\ref{eqnSdeFourier}). Using $N >
1$ and Young's inequality in~(\ref{eqndVHatNSq1}) gives
%
%
\begin{eqnarray}\label{inequal}
\del_t \vert\hat{v}_t (n)\vert^{2} &\leq&- 2 \pi^{2} n^{2}
\vert\hat{v}_t (n)\vert^{2} + 2 \pi n \vert\hat{v}_t (n)\vert
\vert B_t(n)\vert\nonumber\\[-8pt]\\[-8pt]
&\leq&-c n^2 \vert\hat{v}_t (n)\vert^{2} + C \vert
B_t(n)\vert^2\nonumber
\end{eqnarray}
almost surely, where, as before $c, C$ are absolute constants
(independent of $u_0, T$), which may change from line to line. Thus,
for any $t_0' \geq0$ we have
%
%
\begin{equation}\label{eqnHatVn}
\vert\hat{v}_t (n)\vert^{2} \leq\vert\hat v_{t_0'}(n)\vert^2 e^{-
n^2 c t} +
C \int_{t_0'}^t e^{-c n^2 (t - t')} \vert B_{t'}(n)\vert^2 \,dt'
\end{equation}
almost surely, by Gronwall's inequality.

By Parseval's identity we know $\vert B_{t}(n)\vert\leq\Vert
v_{t}\Vert_{L^2}^{2}$, and by conservation of energy [equation (\ref
{eqnL2normV})] this gives $\vert B_t(n)\vert\leq\Vert u_{0}\Vert_{L^2}^{2}$
almost surely. Thus the second term in the previous inequality is
bounded from above by $\frac{C}{cn^2} \Vert u_0\Vert_{L^2}^4$. Since
$\vert\hat u_{t_0'}\vert^2 \leq\Vert u_{t_0'}\Vert_{L^2}^2 \leq
\Vert u_0\Vert_{L^2}^2$,
given a lower bound on $t - t_0'$, we can certainly arrange the same
inequality for the first term. Thus choosing $t_1 = \frac{T}{2}$, for
instance, and applying~(\ref{eqnHatVn}) with $t_0' = 0$, we obtain
%
%
\begin{equation}\label{eqn-2}
{\sup_{t \geq t_1}} \vert\hat{v}_t (n)\vert^{2} \leq\frac{C_0}{n^{2}}
\end{equation}
almost surely, where $C_0 = C_0( \Vert u_0\Vert_{L^2}, T )$ is some constant.

Now we bootstrap, and use~(\ref{eqn-2}) to obtain a better estimate on
$B_t$. Assume inductively that for some $\alpha\in\Z^+$, and
$t_{\alpha} = \frac{\alpha}{\alpha+1} T$, we have
%
%
\begin{equation}\label{eqnInductiveVhatAssumption}
{\sup_{t \geq t_{\alpha}} }\vert\hat{v}_t (n)\vert^{2} \leq\frac
{C_\alpha
}{\vert n\vert^{\alpha+ 1}}
\end{equation}
almost surely. Here $C_\alpha= C_\alpha(\Vert u_0\Vert_{L^2}, T,
\alpha)$
is a constant which we allow to change from line to line if necessary.
We will now establish~(\ref{eqnInductiveVhatAssumption}) for $\alpha+
1$. Note that almost surely, for any $t > t_\alpha$, we have
%
%
\begin{eqnarray}\label{eqn-3}\qquad
\vert B_{t}(n)\vert
&\leq&\sum_{m \in\Z} \vert\hat v_t(n-m)\vert\vert\hat
v_t(m)\vert
\leq2 \sum_{\vert m\vert\geq\vert n\vert/2} \vert\hat
v_t(n-m)\vert\vert\hat v_t(m)\vert\nonumber\\
&\leq&2 \Vert v_{t}\Vert_{L^2} \biggl( \sum_{\vert m\vert\geq\vert
n\vert/2}
\vert\hat v_t(m)\vert^{2}\biggr)^{1/2}
\leq2 \Vert u_{0}\Vert_{L^2} \biggl( \sum_{\vert m\vert\geq\vert
n\vert/2}
\frac
{C_\alpha}{m^{\alpha}} \biggr)^{1/2}\\
&\leq&\frac{C_\alpha}{\vert n\vert^{({\alpha-1})/{2}} }.\nonumber
\end{eqnarray}
Now returning to~(\ref{eqnHatVn}) and choosing $t_0' = t_\alpha$, we
see that the second term is bounded by $\frac{C}{cn^2}\frac
{C_0}{n^{\alpha- 1}} = \frac{C C_0}{c n^{\alpha+ 1}}$. For any $t
\geq t_{\alpha+1}$, we can certainly arrange the same inequality for
the first term, and hence this establishes (\ref
{eqnInductiveVhatAssumption}) for $\alpha+ 1$.

Finally note that if~(\ref{eqnInductiveVhatAssumption}) holds for
$\alpha$, then~(\ref{eqnLongTimeHsBoundV}) holds for any $s < \frac
{\alpha}{2} - 1$, completing the proof.
\end{pf*}
\section{\texorpdfstring{Proof of Lemma \protect\ref{lmaCnBound}}{Proof of Lemma 2.7}}\label{sxnCnBound}
In this section, we prove the almost sure $C^n(\T)$ bounds on $u$
stated in Lemma~\ref{lmaCnBound}. We need a few preliminary results first.
\begin{proposition}[(Local existence without resetting)]\label
{ppnLocalExistence}
Let $u_{t_0}$ be a $C^1(\T)$ valued $\mathcal F_{t_0}$-mea\-sur\-able
random variable such that
\[
\Vert u_{t_0}\Vert_{C^{1}} \leq U_1^0
\]
almost surely. There exists $T_0 = T_0( U_1^0 )$, independent of $N$,
such that the solution to~(\ref{eqndXiNt0})--(\ref{eqnUnt0}) exists on
the interval
$[t_0, t_0 + T_0]$. Further if for some $n \geq1$, $u_{t_0}$ is a
$C^n(\T)$ valued, $\mathcal F_{t_0}$-measurable random variable with
\[
\Vert u_{t_0}\Vert_{C^{n}} \leq U_n^0
\]
almost surely, then there exists $U_n = U_n( U_n^0, n )$ such that
%
%
\begin{equation}\label{eqnUCnbound}
{\sup_{t_0 \leq t \leq t_0 + T_0}} \Vert u_t\Vert_{C^{n}} \leq U_n
\end{equation}
almost surely.
\end{proposition}
%

Proposition~\ref{ppnLocalExistence} can be proved using a standard
Picard iteration. A proof of the analogous result for the
Navier--Stokes equations appeared in the Appendix of \cite
{bblParticleMethod} (see also~\cite{bblSPerturb,bblThesis}). The proof
of~\ref{ppnLocalExistence} is very similar, and we do not provide it here.
\begin{lemma}\label{lemmaaa}
Let $I\dvtx\R\to\R$ denote the identity function, $d \in[0 , 1)$ and let
$\lambda\in C^n(\T)$ be a periodic function such that $\Vert\delx
\lambda\Vert_{L^\infty} \leq d$. Then there exists a constant
$c_{n-1} = c_{n-1}
(\Vert\delx^{n-1} \lambda\Vert_{L^\infty}, d, n)$ such that for
any $f \in
C^n(\R)$,
%
%
\begin{eqnarray}\quad
\label{eqnDnfoIplusLambda} \Vert\delx^n [f\circ(I + \lambda
)] \Vert_{L^\infty} &\leq&\Vert\delx^n f\Vert_{L^\infty}
( 1 + \Vert\delx\lambda\Vert_{L^\infty} )^n +
c_{n-1} \Vert\delx^n \lambda\Vert_{L^\infty},\\
\label{eqnDnIplusLambdaInv} \Vert\delx^n (I + \lambda)\inv\Vert
_{L^\infty}
&\leq&
c_{n-1} \Vert\delx^n \lambda\Vert_{L^\infty},
\end{eqnarray}
for $n > 1$.
\end{lemma}
%
%
\begin{remark*}
Note that since $\Vert\delx\lambda\Vert_{L^\infty} < 1$, the
function $I +
\lambda$ is a $C^1(\R)$ diffeomorphism of $\R$. The notation $(I +
\lambda)\inv$ in~(\ref{eqnDnIplusLambdaInv}) refers to the inverse of
the $C^1(\R)$ diffeomorphism $I + \lambda$.
\end{remark*}
\begin{pf*}{Proof of Lemma~\ref{lemmaaa}}
First note that we can view $\lambda$ as a periodic function (with
period $1$) in $C^n(\R)$. Further, by the mean value theorem, for any
$k \geq1$ there exists $x \in\T$ such that $\delx^k \lambda(x) = 0$.
Thus for any $k \in\{1,\ldots, n\}$, we have $ \vert\delx^k \lambda
\vert
\leq c(n)\Vert\delx^n \lambda\Vert_{L^\infty}$, for some constant $c(n)$
depending only on $n$. (For $k = 0$, we need to subtract the mean of
$\lambda$ for this bound to be valid.)

Now for any two $f, g \in C^n(\R)$, we have
%
%
\begin{equation}\label{eqnDnfog}
\delx^n (f \circ g) = \sum_{m=1}^n (\delx^m f) \circ g
\mathop{\sum
_{k_1 + \cdots+ k_m = n}}_{k_i \geq1} \prod_{i=1}^m \delx
^{k_i} g.
\end{equation}
To prove~(\ref{eqnDnfoIplusLambda}), we set $g = I + \lambda$. The term
in~(\ref{eqnDnfog}) corresponding to $m = n$ gives the first term
of~\ref{eqnDnfoIplusLambda}. When $m < n$, we notice that $k_i > 1$ for
at least one $i$, and $k_j \leq n-1$ for all other $j$. Thus
$\Vert\delx^{k_i} (I + \lambda)\Vert_{L^\infty} = \Vert\delx
^{k_i} \lambda\Vert_{L^\infty} \leq c(n) \Vert\delx^n \lambda
\Vert_{L^\infty}$. The remaining terms $\delx
^{k_j} (I + \lambda)$, $j \neq i$ in the product can be bounded by
$c_{n-1}$. This proves~(\ref{eqnDnfoIplusLambda}).

For~(\ref{eqnDnIplusLambdaInv}), set $X = I + \lambda$ and $A = X\inv$.
Since $n > 1$, $\delx^n (A \circ X) \equiv0$, and using (\ref
{eqnDnfog}) we obtain
\[
\delx^n A \at_X = \frac{-1}{(\delx X)^n} \sum_{m = 1}^{n-1} \delx
^m A
\at_X \mathop{\sum_{k_1 + \cdots+ k_m = n}}_{k_i \geq1} \prod
_{i=1}^m \delx^{k_i} A.
\]
By induction, one can assume that $\Vert\delx^m A\Vert_{L^\infty}
\leq c_{n-1}$
for all $m \leq n-1$. Since $d < 1$, $\frac{1}{\Vert\delx X\Vert
_{L^\infty}}
\leq
\frac{1}{1-d}$, and remaining terms can be bounded by the same argument
as before. This proves~(\ref{eqnDnIplusLambdaInv}).
\end{pf*}
%
%
\begin{lemma}\label{lmaCnGrowthRate}
Let $n \in\N$, $u_{t_0}$ be a bounded, $C^n(\T)$ valued, $\mathcal
F_{t_0}$-measurable random variable. For $k \in\{0,\ldots, k\}$, let
$U^0_k$ be a constant such that $\Vert u_{t_0}\Vert_{C^{k}} \leq
U^0_k$ almost
surely. Let $u$ be the solution of~(\ref{eqndXiNt0})--(\ref
{eqnUnt0}) with initial data
$u_t = u_{t_0}$ when $t = t_0$. If $n > 1$, there exists $\Omega' \in
\mathcal F_{t_0}$ with $P(\Omega') = 1$, $T_0 = T_0(U^0_1) > t_0$ and a
constant $c_{n-1} = c_{n-1}( U^0_{n-1}, n)$ such that
%
%
\begin{equation}\label{eqnDnutbound}
\Vert\delx^n u_t(\omega') \Vert_{L^\infty} \leq U^0_n \bigl( 1 +
c_{n-1} (t -
t_0) \bigr)
\end{equation}
for all $\omega' \in\Omega'$, $t \in[t_0, t_0 + T_0]$. For $n =
1$,~(\ref{eqnDnutbound}) holds with $c_0$ to be an absolute constant.
\end{lemma}
\begin{pf}
For simplicity, we assume $t_0 = 0$. One can check that this assumption
does not affect our proof below. Our first step is to obtain almost
sure $C^1(\T)$ estimates on the Eulerian and Lagrangian displacements.
Throughout this section, we use the convention that $c_{n-1} =
c_{n-1}(U^0_{n-1}, n)$ is a constant depending only on $n$ and
$U^0_{n-1}$ (or an absolute constant for $n=1$), which can change from
line to line.

Let $T_0 = T_0( U^0_1 )$ be the local existence time given by
Proposition~\ref{ppnLocalExistence}, and $c_1 = c_1(U^0_1)$ the almost
sure bound on $\Vert u_t\Vert_{C^{1}}$ from~(\ref{eqnUCnbound}). Let
$I\dvtx\R
\to\R$ be the identity map, $X^i$, $A^i$, respectively, be as in (\ref
{eqndXiNt0}),~(\ref{eqnAnt0}), with $\tau= T_0$. Define $\lambda^i_t =
X^i_t - I$, $\ell^i_t = A^i_t - I$.

Differentiating~(\ref{eqndXiNt0}) with respect to $x$ we obtain
\[
\Vert\delx\lambda^i_t \Vert_{L^\infty} \leq\int_0^t \Vert\delx
u_s\Vert_{L^\infty}
( 1 + \Vert\delx\lambda^i_t \Vert_{L^\infty})
\]
almost surely, for $t \in[0, T_0]$. By Gronwall's lemma,
\[
\Vert\delx\lambda^i_t \Vert_{L^\infty} \leq e^{c_1 t} \int_0^t
\Vert\delx u_s\Vert_{L^\infty} \,ds \qquad\as
\]
for $t \in[0, T_0]$. Recall $t \leq T_0$, $c_1$ only depends on
$U^1_0$, and for all $s \leq T_0$, $\Vert\delx u_s\Vert_{L^\infty}
\leq U^0_1$
almost surely. Thus, as $T_0$ is allowed to depend on $U^0_1$, by
making $T_0$ smaller if necessary we can arrange
%
%
\begin{equation}\label{eqnD1lambda}
\Vert\delx\lambda^i_t \Vert_{L^\infty} \leq c_0 \int_0^t \Vert
\delx u_s\Vert_{L^\infty}
\,ds
\quad\mbox{and}\quad
\sup_{0 \leq t \leq T_0} \Vert\delx\lambda^i_t \Vert_{L^\infty}
\leq\frac{1}{2}
\end{equation}
almost surely, for some absolute constant $c_0$. Now
\begin{eqnarray*}
\delx\ell^i_t &=& \delx A^i_t -1 = \frac{1}{ (\delx X^i_t)
\circ A^i_t} - 1 \\
&=& - \frac{ (\delx\lambda^i_t) \circ A^i_t
}{ 1 + (\delx\lambda^i_t) \circ A^i_t }
\end{eqnarray*}
almost surely. Thus we must have
%
%
\begin{equation}\label{eqnD1ell}
\Vert\delx\ell^i_t\Vert_{L^\infty} \leq2 \Vert\delx\lambda
^i_t\Vert_{L^\infty}
\end{equation}
almost surely for $t \in[0, T_0]$. Using~(\ref{eqnUnt0}) and (\ref
{eqnD1ell}) we have
\begin{eqnarray*}
\Vert\delx u_t \Vert_{L^\infty} &\leq&\frac{1}{N} \sum_{i=1}^N
\Vert\delx u_0 \Vert_{L^\infty} ( 1 + \Vert\delx\ell^i_t
\Vert_{L^\infty} )\\
& \leq&\frac{1}{N} \sum_{i=1}^N \Vert\delx u_0 \Vert_{L^\infty}
( 1 + 2
\Vert\delx\lambda^i_t \Vert_{L^\infty} ) \\
& \leq&\Vert\delx u_0 \Vert_{L^\infty} + 2 c_0 \int_0^t \Vert\delx
u_s\Vert_{L^\infty}
\,ds
\end{eqnarray*}
almost surely for $t \in[0, T_0]$. This proves~(\ref{eqnDnutbound})
for $n=1$.

For $n > 1$, local existence (Proposition~\ref{ppnLocalExistence})
guarantees that $\Vert u_t \Vert_{C^{n-1}} \leq c_{n-1}$ almost surely
for $t
\in[0, T_0]$, where $c_{n-1} = c_{n-1}( U^0_{n-1}, n)$. Assume by
induction that the bound~(\ref{eqnDnutbound}) holds for some integer
$n-1$. This bound and equation~(\ref{eqndXiNt0}) immediately imply that
$\Vert\delx\lambda^i_t \Vert_{C^{n-2}} \leq c_{n-1}$ almost surely%
\footnote{We remark that our somewhat unusual notation $\Vert\delx
\lambda^i_t\Vert_{C^{n-2}}$ instead of $\Vert\lambda^i_t\Vert
_{C^{n-1}}$ is necessary.
This is because it is impossible to obtain almost sure bounds on
$\Vert\lambda^i_t\Vert_{L^\infty}$. However, as our argument shows,
we can obtain
almost sure bounds on $\Vert\delx^k \lambda^i_t\Vert_{L^\infty}$
for any $k
\geq1$.}
for $t \in[0, T_0]$. Equations~(\ref{eqnD1lambda}) and (\ref
{eqnDnIplusLambdaInv}) will imply $\Vert\delx\ell^i_t\Vert
_{C^{n-2}} \leq
c_{n-1}$ almost surely for $t \in[0, T_0]$.

Thus using equations~(\ref{eqndXiNt0}) and~(\ref{eqnDnfoIplusLambda})
we obtain
\begin{eqnarray*}
\Vert\delx^n \lambda^i_t \Vert_{L^\infty} &\leq&\int_0^t \Vert
\delx^n [ u_s \circ(I + \lambda^i_s)]\Vert_{L^\infty
} \,ds\\
&\leq& c_{n-1} \int_0^t [ \Vert\delx^n u_s\Vert_{L^\infty} +
\Vert\delx^n \lambda^i_s\Vert_{L^\infty} ] \,ds
\end{eqnarray*}
almost surely. Using Gronwall's lemma this implies
%
%
\begin{equation}\label{eqnLambdaBound}
\Vert\delx^n \lambda^i_t \Vert_{L^\infty} \leq c_{n-1} \int_0^t
\Vert\delx^n u_s\Vert_{L^\infty} \,ds
\end{equation}
almost surely. Here we absorbed the constant $e^{c_{n-1} t}$ into
$c_{n-1}$, which is valid as $t \leq T_0 = T_0( U^0_1 )$. Now
\begin{eqnarray*}
\Vert\delx^n u_t\Vert_{L^\infty} &\leq&\frac{1}{N} \sum_{i = 1}^N
\bigl(
\Vert\delx^n u_0\Vert_{L^\infty} ( 1 + \Vert\delx\ell
^i_t\Vert_{L^\infty} )^n + c_{n-1}
\Vert\delx^n \ell^i_t\Vert_{L^\infty} \bigr)\\
&\leq&\frac{1}{N} \sum_{i=1}^N \bigl( \Vert\delx^n u_0 \Vert
_{L^\infty} ( 1
+ 2
\Vert\del_x \lambda^i_t \Vert_{L^\infty} )^n + c_{n-1} \Vert\delx
^n \lambda^i_t\Vert_{L^\infty}\bigr)\\
&\leq&\Vert\delx^n u_0 \Vert_{L^\infty} ( 1 + c_{n-1} t ) + c_{n-1}
\int
_0^t\Vert\delx^n u_s\Vert_{L^\infty} \,ds
\end{eqnarray*}
almost surely, where we used~(\ref{eqnDnIplusLambdaInv}) and (\ref
{eqnD1ell}) to obtain the second inequality, and equations (\ref
{eqnD1lambda}) and~(\ref{eqnLambdaBound}) to obtain the third
inequality. Now Gronwall's lemma gives~(\ref{eqnDnutbound}), where we
again absorb the exponential factor $e^{c_{n-1} t}$ into $(1 + c_{n-1}
t)$, by replacing $c_{n-1}$ with a larger constant, which by our
convention we still denote by $c_{n-1}$.
\end{pf}
%
%
\begin{pf*}{Proof of Lemma~\ref{lmaCnBound}}
By Proposition~\ref{ppnLocalExistence}, existence will follow if we
establish~(\ref{eqnUdeltatCn}) for $n=1$. We prove~(\ref{eqnUdeltatCn})
by induction. Since the constant $c_0$ in Lemma~\ref{lmaCnGrowthRate}
is absolute, the proof for $n = 1$ is identical to the proof of the
inductive step. Thus we only prove the inductive step.

Assume that~(\ref{eqnUdeltatCn}) holds for $n-1$, choose $c_{n-1} =
c_{n-1}(U_{n-1})$ to be the constant from Lemma~\ref{lmaCnGrowthRate}.
Thus whenever $\delta_t < T_0$,
%
%
\begin{equation}\label{eqnCknormUk1}
\bigl\Vert\delx^n u^{\delta_t}_{(k+1)\delta_t}\bigr\Vert_{L^\infty} \leq(1
+ c_{n-1}
\delta
_t) \Vert\delx^n u^{\delta_t}_{k\delta_t}\Vert_{L^\infty}
\qquad\as
\end{equation}
holds for all $k \leq\frac{T_0}{\delta_t}$. Iterating this we have
\[
\Vert\delx^n u^{\delta_t}_t\Vert_{L^\infty} \leq(1 +
c_{n-1} \delta
_t
)^{T_0/\delta_t} \Vert\delx^n u_0\Vert_{L^\infty} \qquad\as
\]
for all $t \leq T_0$. Thus we choose $U_n$ to be given by
\[
U_n = \Vert\delx^n u_0\Vert_{L^\infty} \sup_{\delta> 0} ( 1
+ c_{n-1}
\delta
)^{T_0/\delta}.
\]
From~(\ref{eqnSPDEu}) we see that $\int_x u^{\delta_t}_t$ is conserved
almost surely. Since $u^{\delta_t}_t$ is periodic, a~bound on
$\Vert\delx^n u^{\delta_t}_t \Vert_{L^\infty}$ will give us a
bound on $\Vert u^{\delta_t}_t\Vert_{C^{n}}$, completing the proof.
\end{pf*}
\section{\texorpdfstring{Proof of Proposition \protect\ref{ppnNToInfinity}}{Proof of Proposition 2.8}}
\label{sxnNToInfinity}

In this section we prove Proposition~\ref{ppnNToInfinity}. We
reintroduce an $N$ as a superscript to explicitly keep track of the
dependence of our processes on $N$, and prove convergence as $N \to
\infty$.
\begin{pf*}{Proof of Proposition~\ref{ppnNToInfinity}}
Let $w_{t}^{N}=v^{N}_{t}-u^b_{t}$. Then~(\ref{eqnViscousBurgers})
and~(\ref{eqnSPDEv}) give
%
%
\begin{equation}\label{eqnSPDEv2}
dw^{N}_t + w^{N}_t \,\delx v^{N}_t \,dt + u^b_t \,\delx w^{N}_t \,dt -
\frac{1}{2} \,\delx^2 w^{N}_t \,dt + \frac{\delx v^{N}_t}{N} \sum
_{j=1}^N dW^j_t = 0.\hspace*{-28pt}
\end{equation}
Thus, by It\^o's formula
\begin{eqnarray*}
&&\frac{1}{2} d \Vert w^{N}_t\Vert^{2}_{L^{2}} + \biggl( \int_\T(
w^{N}_t)^{2} \,\delx v^{N}_t dx \biggr) \,dt\\
&&\quad{} + \biggl(\int_\T u^b_t w^N_t\,
\delx w^N_t dx\biggr) \,dt + \frac{1}{2} \Vert\delx
w^{N}_t\Vert^{2}_{L^{2}} \,dt\\
&&\quad{}
+ \biggl( \int_\T w^{N}_t \,\delx v^{N}_t dx \biggr) \Biggl( \frac{1}{N}
\sum_{j=1}^N dW^j_t \Biggr) \,dt \\
&&\qquad= \frac{1}{2 N} \Vert\delx
v^{N}_t\Vert^{2}_{L^{2}} \,dt.
\end{eqnarray*}
Taking expectations and integrating by parts we obtain
\begin{eqnarray*}
&&\partial_{t} \E\Vert w^{N}_t\Vert^{2}_{L^{2}}
+ \E\biggl[ \int_{\T} ( w^{N}_t)^{2} ( 2 \,\delx v^{N}_t - \delx
u^b_t ) \,dx \biggr]
+ \E\Vert\delx w^{N}_t\Vert^{2}_{L^{2}}
\\
&&\qquad=
\frac{1}{N} \E\Vert\delx v^{N}_t\Vert^{2}_{L^{2}}.
\end{eqnarray*}
By Lemma~\ref{lmaH2boundV} and the Sobolev embedding theorem, there
exists a constant $C = C(s, \Vert u_0\Vert_{H^{s}})$, independent of $N$,
such that
\[
{\sup_{t \geq0}} \Vert\delx v_t\Vert_{L^\infty} \leq C
\]
almost surely. It is well known that the same estimate holds for $\delx
u^b_t$. Further, since $\E\Vert\delx v_t\Vert_{L^2}^2 \leq\sup
_\Omega
\Vert\delx v\Vert_{L^\infty}^2$, making $C$ larger if necessary we have
\[
\sup_{t \geq0} \E\Vert\delx v_t\Vert_{L^2}^2 \leq C.
\]
Thus
\[
\partial_{t} \E\Vert w^{N}_t\Vert^{2}_{L^{2}} \leq V \biggl( \E
\Vert w^{N}_t\Vert^{2}_{L^{2}} + \frac{1}{N}\biggr).
\]
and, since $w_0 = 0$, Gronwall's lemma gives
\[
\E\Vert w^{N}_t\Vert^{2}_{L^{2}} \leq\frac{1}{CN} ( e^{C t} - 1
)
\]
finishing the proof.
\end{pf*}

\section*{Acknowledgments}
The authors would like to thank the referee for his
insightful comments about the first version of this paper.


%
\printaddresses

\end{document}